\documentclass[a4paper,abstract=true]{scrartcl}
\usepackage[utf8]{inputenc}

\usepackage[UKenglish]{babel} 
\usepackage{amsmath, amsthm, amssymb} 
\usepackage{csquotes, geometry, enumitem} 
\usepackage{xspace, xparse, graphicx} 
\usepackage{tabu} 
\usepackage{tikz-cd}
\usepackage{url}
\usepackage{appendix}

\usepackage[hidelinks]{hyperref}

\theoremstyle{definition}
\newtheorem{defn}{Definition}[section]
\newtheorem{rmk}[defn]{Remark}
\newtheorem{ex}[defn]{Example} 

\theoremstyle{plain}
\newtheorem{thm}[defn]{Theorem}
\newtheorem*{thm*}{Theorem}
\newtheorem{lem}[defn]{Lemma}
\newtheorem{prop}[defn]{Proposition}
\newtheorem*{prop*}{Proposition}
\newtheorem{cor}[defn]{Corollary}

\NewDocumentCommand\party{mmO{}}{\frac{\partial^{#3} #1}{\partial^{#3} #2}}

\NewDocumentCommand\del{}{\partial} 
\NewDocumentCommand\delb{}{\bar{\partial}}

\newcommand{\cM}{\mathcal{M}}
\newcommand{\cL}{\mathcal{L}}

\newcommand{\ZZ}{\ensuremath{\mathbb Z}}
\newcommand{\CC}{\ensuremath{\mathbb C}}
\newcommand{\RR}{\ensuremath{\mathbb R}}
\newcommand{\TT}{\ensuremath{\mathbb T}}

\newcommand{\pd}[1]{\frac{\partial}{\partial #1}} 
\newcommand{\graph}{\mathrm{graph}}

\newcommand{\hor}{\ensuremath{\mathrm{hor}}}

\title{Coisotropic branes in symplectic manifolds}
\author{Charlotte Kirchhoff-Lukat\thanks{KU Leuven, Department of Mathematics, Celestijnenlaan 200B box 2400, BE-3001 Leuven, Belgium.
Email: 
\texttt{c.kirchhofflukat@googlemail.com}
}, Marco Zambon\thanks{KU Leuven, Department of Mathematics, Celestijnenlaan 200B box 2400, BE-3001 Leuven, Belgium.
Email: 
\texttt{marco.zambon@kuleuven.be}}}
\date{}

\begin{document}

\maketitle

\begin{abstract} 
A brane in a symplectic manifold is a coisotropic submanifold $Y$ endowed with a compatible closed 2-form $F$, which together induce a transverse complex structure.  
For a specific class of branes we give  an explicit description of branes nearby a given one, and for arbitrary branes we describe the infinitesimal deformations and provide an associated cochain complex.
As an application, we determine to what extent coisotropic submanifolds near a given brane admit brane structures.
\end{abstract}

\setcounter{tocdepth}{2}
\tableofcontents

\section{Introduction}

This paper is 
devoted to the geometry of branes, as introduced in \cite{GualtieriThesis},
in the setting of symplectic geometry. Given an ambient symplectic manifold $(M,\omega_M)$, a \emph{brane} consists of a coisotropic submanifold $Y$ together with a closed $2$-form $F$ on $Y$ such that
\begin{itemize}
    \item[a)] $F$ and $\omega:=\iota^*_Y\omega_M$ have the same (constant rank) kernel $E$, 
    \item[b)]  the 2-forms $F$ and $\omega$ induce a complex structure transversely to $E$.
\end{itemize}
(See Definition \ref{def:brane} for the precise statement.) When the quotient of $Y$ by the involutive distribution $E$ is smooth, it inherits a complex structure and a holomorphic symplectic form. In the presence of a brane, there is thus an interaction between the symplectic geometry of the ambient manifold and complex geometry. 
 
A brane is  a coisotropic submanifold $Y$ together with a  piece of geometric data, namely a 2-form $F$; such a 2-form, when it exists, is far from being unique. Given a brane, we are interested in the space of nearby branes, for instance to produce new examples. The deformation theory of coisotropic submanifolds   was studied in \cite{OP} and is non-linear, thus much more involved than for the special case of Lagrangian submanifolds. For the deformation theory of branes, we also have to keep track of the geometric data given by the 2-forms. The deformation theory of branes is interesting even when $Y=M$ (\emph{space-filling branes}), i.e.\ when only the 2-form $F$ varies. 

{Branes arise naturally in generalized complex geometry \cite{GualtieriThesis}.} 
We remark that there is also a notion of prequantized brane \cite[\S 4]{KO}, which involves a choice of line bundle with connection. {It was originally motivated by Mirror Symmetry, and is relevant for the quantization of symplectic manifolds \cite{GukovWitten}}. We will not address prequantized branes in this paper, {but will do so in a forthcoming work}.

\bigskip
\noindent{\textbf{Main results}}

Let $(Y,F)$ be a brane in a symplectic manifold. The first aim of this paper is to understand the space of nearby branes. We obtain some partial results in this direction:

\begin{itemize}
    \item For a specific class of codimension $1$ branes of the form $Y=N\times S^1$, where $N$ is a symplectic manifold, we give a description of the nearby branes in terms of the symplectic geometry of $N$ alone.
See \S \ref{subsec:nearby}-\ref{subsec:exobProp}, in particular: 

\begin{prop*}[Proposition \ref{prop:timeone}]
 Fix a symplectic manifold
$(N,\omega_N)$, and  consider the coisotropic submanifold $Y:=N\times S^1\times\{0\}$ of $(N\times S^1\times\RR,\omega_N\times \omega_{T^*S^1})$. Let $f\in C^{\infty}_c(Y)$.
There is a bijection between 
\begin{itemize}
\item [a)] brane structures on   $\graph(f)$, and
 \item [b)] 
space-filling brane structures on $N$ which are
 preserved by $\Phi^1_f$.
\end{itemize} 
Here $\Phi^1_f$ is a certain   symplectomorphism of $(N,\omega_N)$ determined by $f$.
\end{prop*}

\item For arbitrary branes $Y$, in \S \ref{subsec:Diracapproach} we  give a description of the infinitesimal deformations, as elements in $\Gamma(E^*)\oplus \Omega^2(Y)$
satisfying certain conditions. 
The description simplifies when an integrability assumption is made, see  Corollary \ref{rem:infdefcodim1}.
For space-filling branes, the infinitesimal deformations are simply the closed real $(1,1)$-forms.

We also display cochain complexes governing the infinitesimal deformations of branes and their hamiltonian symmetries, in \S \ref{subsec:cochaincomplexes}. They are particularly simple for space-filling branes, see Proposition  \ref{prop:cplxspace-filling}.
{For space-filling branes we provide an obstruction map for the prolongation of an infinitesimal deformation to a path of branes. We show that on the $4$-torus the obstruction map vanishes and  all infinitesimal deformations can be prolonged, thus providing a construction of space-filling branes on the $4$-torus.}
\end{itemize}

To the best of our knowledge, the literature so far considered only \emph{infinitesimal} deformations of branes \cite{KoerMar}\cite{Collier}, while we investigate also \emph{actual} deformations.
Further, the infinitesimal deformations we consider here are more general than those considered in the literature so far, see Remark \ref{complit}.

\bigskip
A second aim of the paper is to determine to what extent the notion of brane is more restrictive than the one of coisotropic submanifold. That is, we want to compare deformations of branes and deformations of coisotropic submanifolds \cite{OP}.

\begin{itemize}
    \item  Given a brane $(Y,F)$ in a symplectic manifold, we ask whether all nearby coisotropic submanifolds admit a brane structure. In other words, we ask whether the forgetful map
 \begin{align*}
   \upsilon\colon \{\text{Branes}\}&\to \{\text{Coisotropic submanifolds}\}\\
 (\widetilde{Y},\widetilde{F})&\mapsto \widetilde{Y}\nonumber
\end{align*}
is surjective near $Y$. 
In general, unsurprisingly, the answer is negative.
This can be seen considering the special case of codimension $1$ branes $N\times S^1$ mentioned above,  see \S \ref{subsec:notbrane}. The proof we provide in Proposition \ref{prop:K3nobrane} however is not straightforward, and indeed relies on  properties of the complex geometry of $K3$ surfaces. It yields:

 \begin{thm*}[Theorem \ref{thm:nobrane}]
There exists a compact brane $Y$ in a symplectic manifold  with this property: 
for every $\epsilon>0$ there is a coisotropic submanifold which is $\epsilon$-close to $Y$ in the $C^2$-sense and  which does not admit a brane structure. 
\end{thm*}

\item  We also ask the infinitesimal version of the above question, in \S \ref{sec:infdefcoiso}. In other words, given a brane $(Y,F)$,  we consider the formal derivative at $(Y,F)$ of the map $\upsilon$, namely a forgetful map
\begin{equation*}
  \Upsilon\colon 
\{\text{Infinitesimal brane deformations}\}  
  \to \{\text{Infinitesimal coisotropic deformations}\}.
\end{equation*}
We ask whether this linear map is surjective in general. We show that the answer is negative, by characterizing explicitly the image of $  \Upsilon$ in the special case of codimension $1$ branes $N\times S^1$ mentioned above, see Proposition \ref{prop:notsurj}.
\end{itemize}
The relation between the two questions addressed in the items above is explained at the beginning of \S \ref{sec:infdefcoiso}:
at least for codimension $1$ branes, the non-surjectivity of $\Upsilon$ implies that there is a \emph{smooth curve} of coisotropics that can not be lifted to a smooth curves of branes.

\bigskip
\noindent\textbf{Notation:} 
Given a 2-form $B\in \Omega^2(Y)$, we abuse notation by denoting by the same symbol the induced bundle map $TY \to T^*Y, v\mapsto \iota_vB$.

\bigskip
\noindent\textbf{Acknowledgements:} 
We thank Marco Gualtieri
for crucial input he provided at the beginning of this project, and Melanie Bertelson, Alberto Cattaneo, Hudson Lima, Karandeep Singh and Konstantin Wernli for useful discussions and observations. {We thank the referees for their constructive comments and suggestions.}

We thank the University of Zurich for hospitality during a research visit in the framework of a UZH Global Strategy and Partnerships Funding Scheme grant.
C.K-L. acknowledges support by the FWO Junior Postdoctoral Fellowship 1263121N and Marie Skłodowska-Curie Fellowship 887857 under the European Union's Horizon 2020 research framework. 
M.Z. acknowledges partial support by
the FWO and FNRS under EOS projects G0H4518N and G0I2222N, by FWO projects G0B3523N and G014726N, and by Methusalem grant METH/21/03 - long term structural funding of the Flemish Government (Belgium).

\section{Branes}\label{sec:branes}

Let $M$ be a symplectic manifold.
We define branes, which are coisotropic submanifolds endowed with an additional piece of data. After spelling out in \S \ref{subsec:space-filling} the  case of space-filling branes -- when the submanifold is the whole of $M$ --, we address some properties and examples of arbitrary branes in \S \ref{subsec:propertyex}. The simple Example \ref{ex:mainex} will be used extensively in  the rest of the paper.

\subsection{Main definitions}

Let $(M,\omega_M)$ be a symplectic manifold.

\begin{defn}
A submanifold $Y$ of $M$ is \emph{coisotropic} if $TY^{\omega_M}\subset TY$,     where $TY^{\omega_M}$ denotes the symplectic orthogonal of $TY$.
\end{defn}
 
Coisotropic submanifolds necessarily have dimension $\ge \frac{1}{2}\dim(M)$. 
Examples include:

\begin{itemize}
\item $M$ itself, and its open subsets,
\item codimension $1$ submanifolds,
\item Lagrangian submanifolds.
\end{itemize}

The intrinsic geometric structure that a coisotropic submanifold $Y$ inherits is $\omega:=\iota^*\omega_M$,
 a  closed, constant rank 2-form on $Y$. 
 In the sequel, we use the term \emph{presymplectic form} to denote a closed $2$-form whose kernel has constant rank.

\begin{rmk}[Coisotropic reduction]
\label{rem:coisored}
Let $Y$ be   a coisotropic submanifold of $(M,\omega_M)$.
We call $TY^{\omega_M}=\ker(\omega)$ the \emph{characteristic distribution}. It  is involutive, hence tangent to a foliation on $Y$, by the Frobenius theorem.
Denote by $\underline{Y}$ the leaf space of the foliation, and assume that it admits a smooth manifold structure such that the quotient map $Y\to \underline{Y}$ is a submersion. Then  $\omega$ descends to symplectic form on $\underline{Y}$. 
In other words, there is a (unique) symplectic form $\underline{\omega}\in \Omega^2(\underline{Y})$ whose pullback to $Y$ is $\omega$.
\begin{center}
\begin{tikzcd} 
 (M,\omega_M) & 	 	 \\ 
(Y,\omega)\arrow[hookrightarrow]{u}{\iota_Y}\arrow{r} & (\underline{Y},\underline{\omega})\\
\end{tikzcd} 
\end{center}
\end{rmk}

If $Y$ is a coisotropic submanifold, we saw above that it inherits the structure of a presymplectic manifold. All presymplectic manifolds arise like this, by Gotay's theorem, which generalizes Weinstein's Lagrangian tubular neighborhood theorem:

\begin{rmk}[Gotay's theorem  \cite{Gotay}]
\label{rem:gotay}
Let $(Y,\omega)$ be a presymplectic manifold. { Denote  by $E:=
ker(\omega)$ its  characteristic distribution.}
Make  a choice of distribution $G$ such that $G\oplus E=TY$. Consider the vector bundle $E^* \overset{q}{\rightarrow} Y$, endowed with the 2-form  $\Omega:=q^*\omega+j^*\omega_{T^*Y}$, where $j\colon E^* \cong G^{\circ}\hookrightarrow T^*Y$ is the inclusion. On  a neighborhood $M$ of the zero section, this 2-form is symplectic, and contains $Y$ as a coisotropic submanifold such that the pullback of $\Omega$ is $\omega$. 

Further, any two such symplectic manifolds   are symplectomorphic in a neighborhood of $Y$. (In particular,
 the choice of complement $G$ is immaterial.)
\end{rmk}   

The following definition of brane appeared in \cite[Example 7.8]{GualtieriThesis}.

\begin{defn}\label{def:brane}
A \textit{brane} 
in a symplectic manifold $(M,\omega_M)$
is a pair $(Y,F)$ consisting of a coisotropic submanifold $Y$ with a closed $2$-form $F\in \Omega^2_{\text{cl}}(Y)$, such that 
\begin{itemize}
    \item[a)]$F$ and $\omega:=\iota^*_Y\omega_M$ have the same (constant rank) kernel $E$, 
    \item[b)] on $TY/E$, the endomorphism $I:=\omega^{-1}\circ F$ satisfies $I^2=-Id_{TY/E}$.
\end{itemize}
{A $2$-form $F$ as above is called  \emph{brane structure}.}
\end{defn}

The endomorphism $I$ defines a complex structure transverse to the characteristic distribution $E$, and $F+i\omega$ induces a holomorphic symplectic structure there, see Remark \ref{rem:fuctorialitybranestrans} and Proposition \ref{prop:charbrane} below. In particular, $rank(TY/E)=dim(M)-2 codim(Y)$ is a multiple of $4$, i.e.\ $\dim(Y)-\frac{1}{2}\dim(M)$ is even.

 Before displaying some properties and examples of branes, let us look at the extreme cases.

\begin{ex}[Lagrangian branes]\label{ex:lagr}
Let $Y$ be a lagrangian submanifold of $(M,\omega_M)$, i.e.\ $TY^{\omega_M}=TY$. Then $\omega:=\iota^*_Y\omega_M$
  is the zero 2-form on $Y$, so $E=TY$. Hence $Y$, together with the   zero 2-form $F=0$, is a brane.
\end{ex}

\subsection{Space-filling   branes}
 \label{subsec:space-filling}

Let $(M,\omega)$ be a symplectic manifold.
A space-filling brane  is a brane of the form $(M,F)$, i.e.\ one supported on the whole of $M$.  Making explicit Definition \ref{def:brane} in this case: A \emph{space-filling brane structure}  consists of a symplectic form $F$  on $M$ such that $I:=\omega^{-1}\circ F$ satisfies $I^2=-Id$.

We recall from  \cite[Prop. 2.6]{KLZTorelli} (see also \cite[\S 4]{KO}):

\begin{prop}\label{prop:charbrane}
Let $(M,\omega)$ be a symplectic manifold. The following structures are in bijection with each other:
\begin{itemize}
    \item[i)] A space-filling brane structure $F$,
    \item[ii)] A complex structure $I$ such that $\omega \circ I$ is skew-symmetric,
    \item[iii)] A complex structure $I$ and a 2-form $F$ such that $F+i\omega$ is holomorphic symplectic w.r.t. $I$.
\end{itemize}
\end{prop}

The bijection between $i)$ and $ii)$ is given by $$F\mapsto I:=\omega^{-1} \circ F \quad \text{ and }\quad I\mapsto F:=\omega \circ I.$$

\begin{ex}
\label{rem:IdfcoordsA}
On $M=\RR^4\cong \CC^2$ take the standard 
complex coordinates $z_1=x_1+i y_1, z_2=x_2+iy_2$. Clearly $dz_1\wedge dz_2$ 
is a  holomorphic symplectic form. Its
real and imaginary parts
are given respectively by
\begin{align*}
F&=dx_1\wedge dx_2-dy_1\wedge dy_2\\
\omega&=dx_1\wedge dy_2+dy_1\wedge dx_2.
\end{align*}
Hence $F$ is a space-filling brane structure on $(\RR^4,\omega)$, by Proposition \ref{prop:charbrane}. Notice that
$I:=\omega^{-1}\circ F$ is the standard complex structure. This example applies also to the 4-torus $\RR^4/\ZZ^4$. 
\end{ex}

\begin{ex}\label{ex:K3}
We call \emph{K3 manifold}
the 4-dimensional, simply connected, compact manifold 
$$M=\{z_0^4+z_1^4+z_2^4+z_3^4=0\}\subset\CC\mathbb{P}^3.$$
 A \emph{complex K3 surface}
is given by $M$  endowed with \emph{any}   complex structure  (we tacitly assume that the complex structure induces the standard orientation on $M$). Any such complex structure supports a holomorphic symplectic form $F+i\omega$,  unique up to non-zero constant multiples.  Hence $(M,\omega)$ is a symplectic manifold for which $F$ is a space-filling brane structure.
\end{ex}

{As pointed out in \cite[\S 2.1]{GukovWitten}, if one  starts from a Hyperk\"ahler manifold $(M,g,I,J,K)$ with corresponding K\"ahler 2-forms $\omega_I,\omega_J,\omega_K$, then $\omega_J+i\omega_K$ is a holomorphic symplectic structure w.r.t. $I$, and  $\omega_J$ is a space-filling brane structure on the symplectic manifold $(M, \omega_K)$ (as one can also see from Proposition \ref{prop:charbrane}).}

\subsection{Basic properties and examples of branes}
\label{subsec:propertyex}

We start with a remark on restrictions.

 \begin{rmk}[Space-filling branes on transversals]\label{rem:fuctorialitybranestrans}
Let $(Y,F)$ be a brane in $(M,\omega_M)$.
Let $S$ be a submanifold of $Y$ which is transverse to the characteristic distribution $E$, in the sense that  $T_yS\oplus E_y=T_yY$ at all $y\in S$. Then  $I:=\omega^{-1}\circ F$ defines a complex structure on $S$. Further $S$, together with the pullbacks of $\omega$ and $F$, is a space-filling brane. 
\end{rmk}

\subsubsection*{Quotients}

We first remark that the quotient of a brane by its characteristic distribution is a space-filling brane, and then construct branes by ``inverting'' this procedure.

\begin{rmk}[Space-filling branes on quotients]\label{rem:fuctorialitybranes}
Let $(Y,F)$ be a brane in $(M,\omega_M)$.
Assume that the quotient $\underline{Y}$ of $Y$ by the foliation tangent to $E$ is smooth. 
Then, by coisotropic reduction as in Remark \ref{rem:coisored}, both $\omega$ and $F$ descend to symplectic forms on $\underline{Y}$, and $ \underline{F}$ is a space-filling brane structure on $(\underline{Y}, \underline{\omega})$. In particular, $\underline{F}+i \underline{\omega}$ is a holomorphic symplectic form there, by Proposition \ref{prop:charbrane}.
\end{rmk}

\begin{ex}[Submersions over space-filling branes]\label{ex:pullback}
Let $(N,\omega_N)$ be a symplectic manifold,  let $\pi\colon Y\to N$ be a submersion. Then $$(Y,\omega:=\pi^*\omega_N)$$ is a presymplectic manifold, with   $\ker \omega=\ker(\pi_*)=:E$.
Upon a choice of complement $G$ of $E$ in $TY$, it is a coisotropic submanifold of a neighborhood $M$ of the zero section in the vector bundle $E^* \overset{q}{\rightarrow} Y$, by Gotay's theorem recalled in Remark \ref{rem:gotay}.

Now assume that  $F_N\in \Omega^2(N)$ is a space-filling brane structure on $(N,\omega_N)$.
Then $(Y,\pi^*F_N)$ is a brane in $M$.

In conclusion, if a symplectic manifold $(N,\omega_N)$ admits a space-filling brane structure and $\pi\colon Y\to N$ is a submersion, then $Y$ admits a brane structure.

For instance, one can take $\pi$ to be any $S^1$ principal bundle over $N$, for instance the trivial bundle $\pi\colon N\times S^1\to N$; if $\omega_N$ represents an integral class, as $\pi$ one can take any prequantization circle bundle of $(N,\omega_N)$.

\begin{center}
\begin{equation}\label{diag:qpi} 
\begin{tikzcd}
 (M,\omega_M)\arrow[d,"q|_M"]& 	 	 \\ 
(Y,\omega)\arrow[r,"\pi"] & (N,\omega_N)\\
\end{tikzcd} 
\end{equation}
\end{center}
\end{ex}

Recall that  any codimension $1$ submanifold is coisotropic.
Not all coisotropic submanifolds carry brane structures, as the next example shows.  

\begin{ex}[The 5-sphere is not a brane]\label{ex:cp2}
Consider the unit sphere $Y=S^5$ in $(\RR^6, \omega_{can})$. Being codimension one, it is a coisotropic submanifold. The quotient $\underline{Y}$ by its characteristic distribution is $(\CC P^2, \underline{\omega})$. If $S^5$ admitted a brane structure, the corresponding  presymplectic form $F$  would descend too, and $\underline{F}+i  \underline{\omega}$ would be a holomorphic symplectic form on $\CC P^2$ (see Remark \ref{rem:fuctorialitybranes}). But $\CC P^2$ does not carry any holomorphic symplectic form.

A direct way to see the previous sentence is as follows. We know that 
$\underline{\omega}$ and $\underline{F}$ are symplectic forms on $\CC P^2$, so their cohomology classes are non-zero. If $\underline{F}+i  \underline{\omega}$ was holomorphic symplectic, they would be linearly independent (if we had $[\underline{F}]=\lambda [\underline{\omega}]$,
then we  would have
 $0= [\underline{\omega}]\wedge [\underline{F}]=\lambda [\underline{\omega}]^2$,
 which is not possible since $\omega$ is a symplectic form on a compact manifold).
But $H^2(  \CC P^2,\RR)$ is one-dimensional, giving a contradiction.
\end{ex} 

\subsubsection*{Products}

 \begin{rmk}[Products]\label{rem:products}
 If $(Y_i,F_i)$ is a brane  in the symplectic manifold $(M_i,\omega_{M_i})$, for $i=1,2$, then the product 
$$(Y_1\times Y_2,F_1\times F_2)$$ is a brane in the symplectic manifold $(M_1\times M_2,\omega_{M_1}\times \omega_{M_2})$. In particular, if $(M_1,F_1)$ is a space-filling brane and $Y_2\subset M_2$ is a Lagrangian submanifold, then $(M_1\times Y_2,F_1\times 0)$ is a brane in the product symplectic manifold. 
 \end{rmk}

We just saw that the product of branes is again a brane. We now use this to see  that \emph{locally} any coisotropic submanifold satisfying the necessary dimensional constraint admits a brane structure.

\begin{rmk}[Local existence]
Let $Y$ be a coisotropic submanifold of  a symplectic manifold $(M,\omega_M)$, and $p\in Y$. There exists a neighborhood $U$ of $p$ in $Y$ such that $(U,\omega:=\iota^*_Y\omega_M)$ carries a brane structure if{f} 
$rank(TY/E)\in 4\ZZ$. 
This dimension condition is necessary, as remarked just after Definition \ref{def:brane}.
The dimension condition is sufficient: taking $U$ to be the domain of a foliated chart for the foliation integrating $E$, and applying Darboux's theorem on $U/E$, we see that there are coordinates $(x_1,\dots,x_{2n},y_1,\dots,y_{2n},t_1,\dots,t_k)$ on $U$ for which $$\omega=\sum
_{j=1}^{n}(dx_{2j-1}\wedge dy_{2j}+dy_{2j-1}\wedge dx_{2j}).$$ A brane structure   on $U$ is then given by  
$$F=\sum_{j=1}^{n}(dx_{2j-1}\wedge dx_{2j}-dy_{2j-1}\wedge dy_{2j}),$$
as one sees using  Example \ref{rem:IdfcoordsA} and Remark \ref{rem:products}.
\end{rmk}

 The following simple example will be important for the rest of this paper.
This example is special case of both Example \ref{ex:pullback} and Remark \ref{rem:products}.

\begin{ex}[Main example]\label{ex:mainex}
Fix a symplectic manifold $(N,\omega_N)$.
Consider the symplectic manifold $T^*S^1=S^1\times \RR$, and denote by $q$ the ``angle coordinate'' on $S^1=\RR/\ZZ$.
We consider
the coisotropic submanifold $$Y:=N\times S^1\times \{0\}$$
in the product symplectic manifold $$(M,\omega_M):=(N\times S^1\times\RR,\omega_N\times \omega_{T^*S^1}).$$

Now let $F_N$ be a space-filling brane structure on $N$. 
Denote by $F$ the pullback of $F_N$ by the projection $Y\to N$.
Then clearly
$$(Y,F)$$
is a brane in $(M,\omega_M)$. Further, all brane structures on $Y$ arise  in this fashion.
\end{ex}

We generalize slightly the above example.
\begin{ex}[Mapping tori and generalizations]\label{ex:mappingtorus}
Fix a symplectic manifold $(N,\omega_N)$ together with a space-filling brane structure $F_N$.

 \begin{enumerate}
\item Let $\phi$ be a symplectomorphism  of $(N,\omega_N)$ preserving  $F_N$.
Consider the mapping torus 
$$Y:=(N\times [0,1])/(x,0)\sim (\phi(x),1),$$
 a fiber bundle over $S^1$.
Since $\phi$ preserves $\omega_N$, the trivial extension  of $\omega_N$ to $N\times [0,1]$ induces a well-defined presymplectic form $\omega$ on $Y$. By Gotay's theorem recalled in Remark \ref{rem:gotay}, 
 $(Y,\omega)$ sits coisotropically in a symplectic manifold $(M,\omega_M)$ (a neighborhood of the zero section of $ Y\times \RR$). 
Further $Y$ is a brane, endowed with the 2-form $F\in \Omega^2(Y)$  induced by the trivial extension of $F_N$ to $N\times [0,1]$.  Notice that the fiber bundle $Y$ comes equipped with a canonical flat connection, whose holonomy is $\phi$.

\item 
Generalizing the above item, let $B$ be a manifold, and fix a point $b_0$. Consider an action $\rho$ of the fundamental group $\pi_1(B,b_0)$ on $(N,\omega_N)$  by symplectomorphisms   preserving  $F_N$. As usual, $\pi_1(B,b_0)$ acts also on the universal cover $\tilde{B}$, by covering transformation. Then
$$Y :=(N\times \tilde{B})/\pi_1(B,b_0),$$ 
  where the quotient is by the diagonal action, 
 is a fiber bundle (with flat connection) over $B$.
 It carries a presymplectic form $\omega$ obtained from $\omega_N$, hence by Gotay's theorem, $(Y,\omega)$ sits coisotropically in a symplectic manifold $(M,\omega_M)$. 
Again, $(Y,F)$ is a brane, where $F\in \Omega^2(Y )$ is induced by the trivial extension of $F_N$ to
$(N\times \tilde{B})$.
\end{enumerate}
\end{ex}

\section{Nearby branes: a codimension one study case}\label{sec:count}

 Given a brane $(Y,F)$ in a symplectic manifold, it is desirable to provide a characterization of all branes supported on submanifolds nearby $Y$.
While at present we are not able to do this in general,
in \S \ref{subsec:nearby}
we do so for a specific class of codimension one branes of the form $Y=N\times S^1$, where $(N,\omega_N)$ is a  symplectic manifold.
We obtain a  description of nearby branes in Proposition \ref{prop:timeone}, thus yielding new examples of branes. We use this to show in  \S \ref{subsec:notbrane}
that coisotropic submanifolds arbitrarily close to a brane $Y$ do not necessarily admit a brane structure, see Theorem \ref{thm:nobrane}.

\subsection{Description of nearby branes}\label{subsec:nearby}

 Fix a symplectic manifold
$(N,\omega_N)$, and consider the coisotropic submanifold $$Y:=N\times S^1\times \{0\}$$ in $(M=Y\times \RR,\omega_M)$, as defined in Example  \ref{ex:mainex}. 
A submanifold $C^1$-close to $Y$ is necessarily the graph of a function $f\in C^{\infty}(Y)$. In this subsection,
we characterize when $\graph(f)$ admits a brane structure, in terms of data on $N$ alone. Notice that we  do not need to assume that $Y$ carries a brane structure.

For any smooth function 
{$f\in C^{\infty}(Y)$,}
the graph is a coisotropic submanifold of $M$.  We want to determine when $\graph(f)$ carries a closed 2-form satisfying the requirements of Definition \ref{def:brane}. 
The presymplectic form $\iota_{\graph(f)}^*\omega_M$, under the natural identification $\graph(f)\cong Y$ given by the projection, becomes \cite[Corollary 3.2]{OP} the presymplectic form $$\omega^f:=\omega_N-d(fdq)$$ on $Y$, where $q$ denotes the standard ``coordinate'' on $S^1{=\RR/\ZZ}$.

\begin{rmk}\label{rem:getbraneonN}
 Notice that the pullback of $\omega^f$ to each slice $N\times \{q\}$ equals $\omega_N$.  Hence, if   $\graph(f)$ admits a brane structure, taking its pullback to $N\times \{0\}$ provides $(N,\omega_N)$ with a space-filling brane structure.
\end{rmk}

The   characteristic distribution of $\omega^f$ is $$E^f:=span\left\{\pd{q}-X^{\omega_N}_f\right\}.$$
Here the vector field $X^{\omega_N}_f$ on $Y$ is determined as follows:
$(X^{\omega_N}_f)|_{N\times \{q\}}:=X^{\omega_N}_{{f_q}}=\omega_N^{-1}df_q$ is the hamiltonian vector field of the restricted function 
$f_q:=f|_{N\times \{q\}}$.

	\begin{figure}[h!]
		\begin{center}
			\begin{tikzpicture}[scale=1]
			\draw[thick,fill=gray!15!white] (-5.5,-1) -- (-3.5,1) -- (4.5,1) --   (2.5,-1) -- cycle; 
			\node[below,right] at (3.8,0) {$Y$};

\draw[purple](-5,-0.5) -- (3,-0.5);	
\draw[purple](-4.5,0) -- (3.5,0);	
\draw[purple](-4,0.5) -- (4,0.5);
\node[below,left] at (0,0.2)[purple] {$E$};
\draw[thick] (-5.5,2) 
.. controls (-5,2.5) and (-4,3.5) .. (-3.5,4);

\draw[thick] (-3.5,4)
.. controls (-1,5.5) and (1,3.5) .. (4.5,4);

\draw[thick] (-5.5,2)
.. controls (-3,3.5) and (-1,1.5) .. (2.5,2);

\draw[thick] (2.5,2)
.. controls (3,2.5) and (4,3.5) .. (4.5,4);
 \node[below,right] at (4,3) {$\graph(f)$};

\draw[thick,cyan] (-3.7,3.8)
.. controls (-3,4.5) and (-1,3) .. (3.5,3);  
             
\draw[thick,cyan] (-4,3.5)
.. controls (-3,4.5) and (-1,2.5) .. (3,2.5);       

\draw[thick,cyan] (-4.5,3)
.. controls (-3,4.5) and (-1,2) .. (2.7,2.2);  
\node[below,left] at (0,3.8)[cyan] {$E^f$};
			\end{tikzpicture}
		\end{center}
	\caption{The submanifolds $Y$ and $\graph(f)$, with the respective characteristic foliations.}
	\end{figure}

 We will make use of the following general fact, whose proof we defer to Appendix \ref{app:proof}.

\begin{lem}\label{rem:Melanie}
 Let $Y$ be a manifold and $F$ a 2-form of constant rank, denote $E:=\ker(F)$. Suppose there is an \emph{involutive} distribution $G$  such that  $E\oplus G=TY$, and denote by $\iota_G\colon G\hookrightarrow TY$ the inclusion.
Then $F$ is closed if{f}:
\begin{itemize}
\item[i)] the distribution $E$ is involutive 
\item[ii)]  the foliated two-form $\iota_{G}^*F$ is invariant w.r.t. the holonomy of $E$, i.e.\: for any $X\in \Gamma(E)$ whose flow preserves\footnote{In other words,  $[X,\Gamma(G)]\subset \Gamma(G)$.} $G$, the flow preserves $\iota_{G}^*F$
\item[iii)] the foliated two-form $\iota_{G}^*F$ is closed (thus symplectic).
\end{itemize}
\end{lem}

\begin{lem}\label{lem:firstbij}

 Fix a symplectic manifold
$(N,\omega_N)$, and as in Example  \ref{ex:mainex} consider the coisotropic submanifold $Y:=N\times S^1\times\{0\}$ of $(N\times S^1\times\RR,\omega_N\times \omega_{T^*S^1})$. Let {$f\in C^{\infty}(Y)$, and assume that the following   diffeomorphism of $N$ exists:} 
\begin{equation}\label{eq:phif}
\Phi^1_f:=\text{   Time $1$ flow of the time-dependent vector field $\{-X^{\omega_N}_{{f_q}}\}_{q\in [0,1]}$}.  
\end{equation}

There is a bijection between 
\begin{itemize}
\item [a)] presymplectic forms $\tilde{F}$ on $Y$ whose kernel is $E^f$, and
 \item [b)] symplectic forms $\tilde{F}_N$ on $N$ which are
 preserved by $\Phi^1_f$. 
\end{itemize}
The bijection reads $\tilde{F}\mapsto \iota_{N\times \{0\}}^*\tilde{F}$. Its inverse maps $\tilde{F}_N$ to the unique extension to a 2-form on $Y$ invariant by the vector field $\pd{q}-X^{\omega_N}_f$.
\end{lem}
Notice that the condition in \emph{b)} is a condition on the preferred  slice $N\times \{0\}$ only.

{If $f$ has compact support, then the time-$1$ flow  $\Phi^1_f$ is guaranteed to exist. Further,  in  \eqref{eq:phif} we
use the identification $S^1=[0,1]/(0\sim 1)$
and thus have $X^{\omega_N}_{{f_0}}=X^{\omega_N}_{{f_1}}$.}  
\begin{proof}
 Suppose $\tilde{F}\in \Omega^2(Y)$ is a 2-form   whose kernel is $E^f$. We apply Lemma \ref{rem:Melanie}, with $G=TN\times \{0\}$. It states that 
 $\tilde{F}$ is closed if{f}
\begin{itemize}
\item[ii)]   the flow of $\pd{q}-X^{\omega_N}_f$ preserves the family of two-forms $\iota_{N\times \{q\}}^*\tilde{F}$
\item[iii)]  the two-forms $\iota_{N\times \{q\}}^*\tilde{F}$ are  closed, for all $q\in S^1$.
\end{itemize}
Notice that in this case, by the   condition ii), $\tilde{F}$ is determined by its pullback
$\iota_{N\times \{0\}}^*\tilde{F}$
 to the preferred slice $N\times \{0\}$. Condition ii) can be rephrased  
 saying that $\iota_{N\times \{0\}}^*\tilde{F}$ is preserved by $\Phi^1_f$, and  condition iii) by saying that $\iota_{N\times \{0\}}^*\tilde{F}$ is closed. Hence $\iota_{N\times \{0\}}^*\tilde{F}$ is a symplectic form on $N$ as in b). This procedure can be inverted. 
 \end{proof}

We can now state the following proposition, which characterizes brane structures on   $\graph(f)$ in terms of data on $N$.

\begin{prop}\label{prop:timeone}
Assume the set-up of Lemma \ref{lem:firstbij}.
There is a bijection between 
\begin{itemize}
\item [a)] brane structures on   $\graph(f)$, and
 \item [b)] 
space-filling brane structures on $N$ which are
 preserved by $\Phi^1_f$.
\end{itemize} 
Under this bijection, from a  $\Phi^1_f$-invariant
space-filling brane structure  
$\tilde{F}_N$ on $N$ we obtain a brane structure on $\graph(f)$ as follows: 
take the unique extension of $\tilde{F}_N$ to a 2-form on $Y$ invariant by the vector field $\pd{q}-X^{\omega_N}_f$; pull back the resulting 2-form via the identification
$\graph(f)\cong Y$.
\end{prop}
\begin{proof}
 Recall that a brane structure on $\graph(f)$  corresponds to a presymplectic form $\tilde{F}$ on $Y$ whose kernel is $E^f$ and  so that on $TY/E^f$ the endomorphism $(\omega^f)^{-1}\circ \tilde{F}$ squares to $-Id_{TY/E^f}$.
 Given such a presymplectic form $\tilde{F}$, apply to it the bijection in Lemma \ref{lem:firstbij}, to obtain
the symplectic form  $\iota_{N\times \{0\}}^*\tilde{F}$ on $N$ 
 preserved by $\Phi^1_f$. It is a space-filling brane structure on $(N,\omega_N)$,  see Remark \ref{rem:getbraneonN}.

 Conversely, given a space-filling brane structure  $\tilde{F}_N$ on $N$ as in b), apply the   bijection of Lemma \ref{lem:firstbij}  to obtain a presymplectic form $\tilde{F}$ on $Y$ whose kernel is $E^f$. It satisfies the quadratic condition in the definition of brane, since  both $\tilde{F}$ and $\omega^f$ are invariant along their (common) kernel $E^f$.
 The description of the bijection follows from Lemma \ref{lem:firstbij}.
\end{proof}

\subsection{Examples and observations for Proposition \ref{prop:timeone}}
\label{subsec:exobProp}

We elaborate on  Proposition \ref{prop:timeone}. This subsection can be safely skipped on a first reading, and is divided in four parts which are independent of each other.

\subsubsection*{Special cases}

We  present some special cases of Proposition \ref{prop:timeone}.
 
\begin{rmk}\label{rem:f=0}
Taking $f=0$,  Proposition \ref{prop:timeone} recovers the fact that brane structures on $Y=N\times S^1\times \{0\}$ are given precisely by the pullbacks to $Y$ of 
space-filling brane structures on $N$.
\end{rmk}

\begin{cor}\label{cor:closed1f}
Assume that  $\tilde{F}_N$ is a space-filling brane structure on $(N,\omega_N)$, denote $\tilde{I}:=\omega^{-1}_N\circ \tilde{F}_N$. Let {$f\in C^{\infty}(Y)$}   so that 
$$ (\tilde{I})^*df_q\in \Omega^1(N) \text{ is closed  for all $q\in [0,1]$},$$
where $f_q:=f|_{N\times \{q\}}$, {and so that the diffeomorphism $\Phi^1_f$ in \eqref{eq:phif} exists.} 
 Then  $\tilde{F}_N$
determines a brane structure  on   $\graph(f)$, via the bijection of Proposition \ref{prop:timeone}.
\end{cor}   
\begin{proof}
Consider condition b) in  Proposition \ref{prop:timeone}.
A special case of that condition occurs when $\tilde{F}_N$ is preserved by the flow at any arbitrary time, i.e.\ $\cL_{X^{\omega_N}_{f_q}}\tilde{F}_N=0$   for all $q$. We have
\begin{equation}\label{eq:LiederXf}
\cL_{X^{\omega_N}_{f_q}}\tilde{F}_N=d \iota_{X^{\omega_N}_{f_q}}\tilde{F}_N=
d (\tilde{I})^*df_q.
\end{equation} 
Here in the second equality  we used that $\iota_{X^{\omega_N}_{f_q}}\tilde{F}_N =(\tilde{I})^*df_q$, 
as a consequence of $\tilde{F}_N\circ \omega^{-1}_N=\tilde{I}^*$
\end{proof}

We present an example for Corollary \ref{cor:closed1f}.
 
\begin{ex}
\label{rem:Idfcoords}
Suppose $N=\RR^4$ (or a 4-torus), with the symplectic structure and space-filling brane as in Example \ref{rem:IdfcoordsA}:
\begin{align*}
F&=dx_1\wedge dx_2-dy_1\wedge dy_2\\
\omega&=dx_1\wedge dy_2+dy_1\wedge dx_2,
\end{align*}
so that $I:=\omega^{-1}\circ F$ is the standard complex structure.
Given a function $f$ on $N$, the condition that    $I^*df$ is closed is quite a stringent one, as it amounts to the four equations
$$\frac{\partial^2 f}{\partial^2 x_1}+\frac{\partial^2 f}{\partial^2 y_1}=0,\; \frac{\partial^2 f}{\partial^2 x_2}+\frac{\partial^2 f}{\partial^2 y_2}=0,\;
\frac{\partial^2 f}{\partial x_1\partial x_2}+\frac{\partial^2 f}{\partial y_1\partial y_2}=0,\;
-\frac{\partial^2 f}{\partial y_1\partial x_2}+\frac{\partial^2 f}{\partial x_1\partial y_2}=0.$$

{Notice that the first two equations imply that $f$ is a harmonic function. Hence on the 4-torus, the only solutions are given by constant functions.
On $\RR^4$, any linear function is a solution. Further, any harmonic function $f(x_1,y_1)$ of $x_1$ and $y_1$ is a solution: the first equation is satisfied by construction, the remaining ones are satisfied since the function does not depend on $x_2$ and $y_2$. (As is well-known, harmonic functions on the $x_1$-$y_1$-plane can
be obtained as the real and imaginary part of holomorphic functions; concrete examples are $x_1^2-y_1^2$, $x_1y_1$, $x_1^4+y_1^4-6x_1^2y_1^2$ and $x_1^3y_1-x_1y_1^3$.)}

{
The  hamiltonian vector field 
of a function  $f(x_1,y_1)$ w.r.t. $\omega$ 
is always complete: we have
 $-X_f^{\omega}=\frac{\partial f}{\partial y_1}\pd{x_2}+\frac{\partial f}{\partial x_1}\pd{y_2}$, so the vector field is tangent 
 to every plane on which $x_2$ and $y_2$ are constant, and restricts to a constant vector field on every such plane. In particular, the time-$1$ flow of $-X_f^{\omega}$ exists. We conclude from Corollary \ref{cor:closed1f} that the graph of any function
on $Y=N\times S^1$ of the form $g(q)f(x_1,y_1)$, where $g$ is a smooth  function on $S^1$ and $f$ a harmonic function on the plane, admits a brane structure. 
}

\end{ex}

\subsubsection*{Two branes which are not symplectomorphic}

 Given a brane $(Y,F)$ in a symplectic manifold, it is clear that applying a symplectomorphism $\phi$ of the ambient symplectic manifold we obtain a new brane, namely $\phi(Y)$ together with the pullback of $F$ by  $(\phi|_Y)^{-1}$.
The next example shows that the branes obtained deforming a given one as in Proposition \ref{prop:timeone} can be genuinely different, in the sense that they   are not related by a symplectomorphism. 

 \begin{ex}\label{ex:differentfoliations}
Assume again the set-up of Example  \ref{ex:mainex}, i.e.\ $Y=N\times S^1\times \{0\}$, specializing to the following
 symplectic manifold 
 (a variation of Example \ref{rem:IdfcoordsA}):
 $$(N,\omega_N):=( S^1\times \RR^3, dx_1\wedge dy_2+dy_1\wedge dx_2),$$ where $x_1$ is the ``angle coordinate'' on $S^1$, and $x_2, y_1, y_2$ are coordinates on $\RR^3$.

For any irrational number $\lambda$, let $f:=\lambda y_2$.
Then
$$\pd{q}-X^{\omega_N}_{f}=\pd{q}-\lambda \pd{x_1}.$$
Notice that the time $1$ flow of this vector field is the translation by $-\lambda \text{ mod }\ZZ$ in the $x_1$-direction. In particular, none of the integral curves of this vector field is periodic, i.e.\ none of them closes up. Therefore $(Y,\omega^f)$ can not be presymplectomorphic to $(Y,\omega)$, since the characteristic leaves of the latter are all circles. (Recall that we introduced $\omega^f$ at the beginning of \S \ref{subsec:nearby}.)

 Now consider the space-filling brane structure $F_N=
 dx_1\wedge dx_2-dy_1\wedge dy_2$
 on $N$. The pullback of $F_N$ to $Y$ by the first projection   makes 
  $(Y,\omega)$ into a brane, see Remark \ref{rem:f=0}.
Further, applying the construction of Proposition \ref{prop:timeone} to $F_N$ also makes $\graph(f)$ into a brane, presymplectomorphic to $(Y,\omega^f)$. (This proposition applies because the time $1$ flow of $-\lambda \pd{x_1}$, which is a translation in the $x_1$ direction in $N$, preserves the constant form $F_N$; alternatively, because $f$ satisfies the PDE's in Example \ref {rem:Idfcoords}.)
 
 In conclusion, $Y$ and $\graph(f)$ are two branes 
 which are not related by a symplectomorphism
of the ambient symplectic manifold. 
\end{ex}

\subsubsection*{All nearby branes are mapping tori}
\begin{prop}\label{prop:mappingtorus}
Fix a symplectic manifold
$(N,\omega_N)$, and 
consider the coisotropic submanifold $Y:=N\times S^1\times\{0\}$ of $(N\times S^1\times\RR,\omega_N\times \omega_{T^*S^1})$. 

All branes $C^1$-close to the coisotropic submanifold $Y$ arise from the mapping torus construction described in Example \ref{ex:mappingtorus} 1).
\end{prop}

\begin{proof}
We use freely the notation introduced at the beginning of \S \ref{subsec:nearby}.
Let $f\in C^{\infty}(Y)$. 
Under the natural identification $\graph(f)\cong Y$, the pullback of the ambient symplectic form corresponds to the presymplectic form $\omega^f$, and  a brane structure on $\graph(f)$  corresponds to a 2-form  $\tilde{F}$ on $Y$ whose kernel is $E^f$ and  so that on $TY/E^f$ the endomorphism $(\omega^f)^{-1}\circ \tilde{F}$ squares to $-Id_{TY/E^f}$. 
Denote by $\tilde{F}_N$ the pullback of $\tilde{F}$ to the   preferred slice $N\times \{0\}$.  Assume that $f$ is sufficiently $C^1$-small, then the time $1$ flow   $\Phi^1_f$ exists, and  preserves 
$\tilde{F}_N$ by
Proposition \ref{prop:timeone}.

 For all $q\in [0,1]$, denote
\begin{equation*} 
\Phi^q_f:=\text{   Time $q$ flow of the time-dependent vector field  $\{-X^{\omega_N}_{{f_q}}\}_{q\in [0,1]}$}.  
\end{equation*}
Consider now the mapping torus $Z$ associated to the diffeomorphism $(\Phi^1_f)^{-1}$ of $N$.
The map
$$\psi\colon Z:= \frac{N\times [0,1]}{(\Phi^1_f (x),0)\sim (x,1)}\to Y=N\times S^1,\;\; (x,q)\mapsto (\Phi^q_f(x),q)$$
is   well-defined and a diffeomorphism. Here we use the identification $S^1=[0,1]/(0\sim 1)$. The map $\psi$ sends the tangent vector $\pd{q}$ at $(x,q)$ to the tangent vector at $\psi(x,q)$ given by $$\frac{d}{ds}|_{s=0}(\Phi^{q+s}_f(x), q+s)=\pd{q}-X^{\omega_N}_{{f_q}}.$$ 
Therefore $\psi$ intertwines
\begin{itemize}
    \item the distribution $\RR\pd{q}$ on $Z$, i.e.\ the kernel
of the 2-forms (induced by) $\pi^*\omega_N$ and $\pi^*\tilde{F}_N$
\end{itemize}
with
\begin{itemize}
\item the distribution $E^f=\RR(\pd{q}-X^{\omega_N}_f)$ on $Y$, i.e.\ the kernel of the 2-forms $\omega^f$ and  $\tilde{F}$.
\end{itemize} 
Both $(\psi^{-1})^* (\pi^*\tilde{F}_N)$ and $\tilde{F}$ are closed 2-forms on $Y=N\times S^1$ with the same kernel, so they are determined by their pullbacks
 to the preferred slice $N\times \{0\}$, by the proof of Lemma \ref{lem:firstbij}. We have $(\psi^{-1})^* (\pi^*\tilde{F}_N)=\tilde{F}$, because
their  pullbacks to $N\times \{0\}$ agree (they both equal $\tilde{F}_N$, as one sees using the fact that $\Phi^1_f$ preserves $\tilde{F}_N$). A similar reasoning shows that $(\psi^{-1})^* (\pi^*\omega_N)= \omega^f$.

This shows that $\psi$ is an isomorphism between  
the presymplectic manifold $(Y,\omega^f)$ together with the 2-form $\tilde{F}_N$, and the 
 presymplectic manifold $(Z, \pi^*\omega_N)$ together with the 2-form  $\pi^*\tilde{F}_N$ (the  mapping torus).
 By the uniqueness in Gotay's theorem, this finishes the proof.
\end{proof}

\subsubsection*{Non-existence of nearby branes}

Roughly speaking, in the setting of Remark \ref{rem:getbraneonN} and Proposition \ref{prop:timeone}, on submanifolds nearby $Y$ there are ``less'' brane structures than on $Y$ itself. This immediately implies a non-existence result for brane structures on nearby submanifolds.
\begin{cor}\label{cor:nobranes}
 Let $(N,\omega_N)$ be a symplectic manifold.
 As in Example  \ref{ex:mainex}, consider the symplectic manifold $(N\times S^1\times\RR,\omega_N\times \omega_{T^*S^1})$ 
 and its  submanifold $Y:=N\times S^1\times\{0\}$. 
 Assume that $(N,\omega_N)$ admits no space-filling brane structure.
 Then, for all $f\in C^{\infty}(Y)$, the submanifold $\graph(f)$ admits no brane structure. 
\end{cor}
An instance of symplectic manifold admitting no space-filling brane structure is $\CC P^2$, as recalled in Example \ref{ex:cp2}.

\begin{proof}
If  $\graph(f)$ admitted a brane structure, by 
Remark \ref{rem:getbraneonN} 
$(N,\omega_N)$ would  admit a space-filling brane structure, contradicting the assumption.
\end{proof}

\begin{rmk}
 In particular, if $Y$ as in Corollary \ref{cor:nobranes} admits no brane structure, then all submanifolds $C^1$-close to $Y$ also do not admit a brane structure. 
 An interpretation of this statement is the following.
Consider the space $\mathcal{C}$ of (necessarily coisotropic) codimension one submanifolds of $M=N\times S^1\times\RR$, endowed with a suitable topology, and the subset $\mathcal{B}$ of codimension one submanifolds  admitting a brane structure. Then $\mathcal{C}\setminus \mathcal{B}$ contains an open subset, because $Y$ has an open neighborhood contained in $\mathcal{C}\setminus \mathcal{B}$.

 The above statement should be compared with Theorem \ref{thm:nobrane} below: for a specific choice of $(N,\omega_N)$, for which $Y=N\times S^1\times\{0\}$ admits a brane structure, the subset $\mathcal{B}$ is not open in $\mathcal{C}$, and indeed  $Y$ is a boundary point of $\mathcal{B}$. 
\end{rmk}

We present a generalization of Corollary \ref{cor:nobranes}, in which $Y$ is not necessarily a product, and which allows to conclude that  all coisotropic submanifolds which are sufficiently $C^1$-close to $Y$ do not admit space-filling brane structures. The set-up is displayed in diagram \eqref{diag:qpi}.

\begin{prop}\label{prop:Ynobranesnearby}
 Let $(N,\omega_N)$ be a symplectic manifold. Let $\pi\colon Y\to N$ be a submersion admitting a section $s$, and denote $E:=\ker \pi_*$.
 Choose a distribution $G$ on $Y$ such that $E\oplus G=TY$ and so that $s(N)$ is tangent to $G$.
 The presymplectic manifold  $(Y,\omega:=\pi^*\omega_N)$ embeds coisotropically in a neighborhood $(M,\omega_M)$ of the zero section of $E^*\to Y$, as in  Example \ref{ex:pullback}, by Gotay's construction.

Suppose that $(N,\omega_N)$ admits no space-filling brane structure. Then, for all sufficiently small sections $r$ of $q\colon M\to Y$, the submanifold $\graph(r)$ of $(M,\omega_M)$ does not admit any brane structure.
\end{prop}
\begin{proof}
Suppose $\graph(r)$ is coisotropic and admits a brane structure.
Since $r$ is small, $r(s(N))$ is a submanifold transverse to the characteristic distribution of $\graph(r)$. By Remark \ref{rem:fuctorialitybranestrans}, $r(s(N))$ admits a space-filling brane structure w.r.t. the pullback of $\omega_M$; equivalently, $s(N)$ admits a space-filling brane structure w.r.t. the pullback of $r^*\omega_M$.

 We know that $r^*\omega_M=\omega-d\bar{r}$ \cite[Corollary 3.2]{OP}, where $\bar{r}\in \Omega^1(Y)$ is the extension of $r  \in \Gamma(E^*)$ annihilating $G$.
 Notice that $\iota_{s(N)}^*\bar{r}=0$, since $s(N)$ is tangent to $G$. Hence $\iota_{s(N)}^*(r^*(\omega_M))=\iota_{s(N)}^*\omega$,
 and under the diffeomorphism $s(N)\cong N$ it corresponds to $\omega_N$. Therefore $(N,\omega_N)$ carries a space-filling brane structure, contradicting our assumption.
\end{proof}

\begin{ex}
i) Let $(N,\omega_N)$ be a symplectic manifold. In Proposition \ref{prop:Ynobranesnearby} one can take $\pi$ to be a vector bundle $\pi\colon V\to N$, and $s$ any section. Notice that in this case $E=\ker(\pi_*)$ is the vertical bundle, which is isomorphic to the pullback vector bundle $\pi^*V$.

For instance, take the symplectic manifold $N=\CC P^2$, which does not admit a space-filling brane structure by Example \ref{ex:cp2}, and let $V$ be the tautological complex line bundle (so $V$ as a manifold is diffeomorphic to the complex blow-up of $\CC^3$ at the origin).
Then no submanifold $C^1$-close to $V$ in $E^*$ admits a brane structure,
by Proposition \ref{prop:Ynobranesnearby}.

ii) The unit sphere $Y=S^5$ in $(\RR^6, \omega_{can})$ does not admit a brane structure, by Example \ref{ex:cp2}. We do not know if the same holds for all submanifolds close enough to $S^5$, since 
Proposition \ref{prop:Ynobranesnearby} does not apply in this case (the projection $S^5\to \CC P^2$ admits no global section).
\end{ex}

\subsection{Coisotropic deformations which are not branes}\label{subsec:notbrane}

In this subsection we exhibit an example of brane admitting small coisotropic deformations which are not branes. In other words, given a brane $(Y,F)$ in a symplectic manifold, we show that in general the forgetful map
 \begin{align}\label{eq:forgetful}
   \upsilon\colon \{\text{Branes}\}&\to \{\text{Coisotropic submanifolds}\}\\
 (\widetilde{Y},\widetilde{F})&\mapsto \widetilde{Y}\nonumber
\end{align}
is not surjective near $Y$. The precise statement is given in Theorem \ref{thm:nobrane}.
 \smallskip

Let $(N,\omega_N)$ be a  symplectic manifold with a space-filling brane $F_N$. We saw in Example  \ref{ex:mainex}  that $Y:=N\times S^1\times \{0\}$, together with the presymplectic form $F$ obtained extending trivially $F_N$ to $Y$, is a brane
in the symplectic manifold $(M,\omega_M):=(N\times S^1\times\RR,\omega_N\times \omega_{T^*S^1})$. 
We apply this to the K3 manifold, following a suggestion of Marco Gualtieri. Recall from Example \ref{ex:K3} that any complex structure on the K3 manifold yields a symplectic form  admitting a space-filling brane.

\begin{prop}\label{prop:K3nobrane}
Let $N$ be the K3 manifold.  Let $\omega_N$ be a symplectic form admitting a space-filling brane $F_N$.
Let $f\in C^{\infty}(Y)$. Then:

$$\text{$\graph(f)$ admits brane structure \quad $\Leftrightarrow$\quad $\Phi^1_f=Id_N$},$$
where $\Phi^1_f$ is the diffeomorphism of $N$ defined in \eqref{eq:phif}.
\end{prop}

\begin{rmk}\label{rmk:K3sympl}
If $\Phi$ is a symplectomorphism defined on a neighborhood of the brane $Y$ in $(M,\omega_M)$, then $\Phi(Y)$ is again a brane, as observed just before Example \ref{ex:differentfoliations}. We do not know if, in the set-up of Proposition \ref{prop:K3nobrane}, all  branes $C^1$-close to $Y$ arise in this way. The infinitesimally version of this statement is true by Remark \ref{rmk:K3const}.
 \end{rmk}

\begin{proof}
If $\Phi^1_f= Id_N$ then $\graph(f)$ admits a brane structure, by Proposition \ref{prop:timeone}.

Conversely, suppose $\graph(f)$ admits  a brane structure. By  Proposition \ref{prop:timeone}, there is a space-filling brane structure $\tilde{F}_N$ on $N$ which is preserved by the  time $1$ flow $\Phi^1_f$ of the time-dependent vector field  $\{-X^{\omega_N}_{{f_q}}\}_{q\in [0,1]}$. Since
this vector field is hamiltonian w.r.t. $\omega_N$, the diffeomorphism $\Phi^1_f$ preserves $\omega_N$. 
 As a consequence, it preserves
the complex structure $\tilde{I}:=\omega^{-1}_N\circ\tilde{F}_N$.

By construction, $\Phi^1_f$   is isotopic to the identity, hence
it induces the identity in cohomology. Proposition \ref{prop:K3id} below implies that 
$\Phi^1_f=Id_N$.  
\end{proof}

\begin{prop}\label{prop:K3id}
\cite[Ch. 15.2, Prop. 2.1]{HuyLecturesK3}
 \cite[Cor. VIII.11.3]{BarthPetersVandeVen1984}
Let $X$ be a complex K3 surface. Let $f$ be an automorphism (i.e.\ it preserves the complex structure). If $f$ induces the identity on $H^2(X, \ZZ)$, then $f = Id_X$. 
\end{prop}

We can now provide a statement about the forgetful map $\upsilon$ in \eqref{eq:forgetful}: 

 \begin{thm}\label{thm:nobrane}
There exists a compact brane $Y$ in a symplectic manifold  with this property: 
for every $\epsilon>0$ there is a coisotropic submanifold which is $\epsilon$-close to $Y$ in the $C^2$-sense and  which does not admit a brane structure. 
\end{thm}
 
\begin{proof}
Take $Y=N\times S^1\times \{0\}$ to be as above, for $N$ the K3 manifold with a symplectic form $\omega_N$ admitting a space-filling brane.
Let $g\in C^{\infty}(N)$ be a non-constant function, and denote by  $f=\pi^*g$ the pullback w.r.t.\ $\pi\colon Y\to N$; then $\Phi^1_f=\Phi^1_g$, where the latter denotes the time one flow of the (time-independent) hamiltonian vector field  $-X^{\omega_N}_g$. 

If $g$ is sufficiently $C^2$-small, the vector field $X^{\omega_N}_g$ is sufficiently $C^1$-small and non-trivial.
Therefore  $\Phi^1_g\neq Id_N$ (this follows for instance from the period bounding lemma by Debord \cite[Proposition 23 in Appendix A]{DebordThesis}). Hence $\graph(f)$ is not a brane, by Proposition \ref{prop:K3nobrane}.
\end{proof}

\section{Infinitesimal deformations of arbitrary branes}\label{sec:infdef}

 In this section we let $(M,\omega_M)$ be a symplectic manifold, and $(Y,F)$ an \emph{arbitrary} brane. In Corollary \ref{cor:infdefdirect} we provide a description of the infinitesimal deformations of $(Y,F)$ as a brane, i.e.\ of the formal tangent space at $(Y,F)$ to the space of branes. We then identify a cochain complex governing these infinitesimal deformations, and whose infinitesimal equivalences of deformations are given by the action of hamiltonian  vector fields.
Finally, we look at the forgetful map from infinitesimal deformations of branes to those of coisotropic submanifolds, and show that it is not surjective in general.

\subsection{Description of infinitesimal deformations}\label{subsec:Diracapproach}

 In this subsection we  consider a whole curve of branes and determine the conditions satisfied by    its initial velocity, in Lemma \ref{lem:infdefdirect}. From this we extrapolate the description of infinitesimal deformations of branes in Corollary \ref{cor:infdefdirect}.
 This description depends on a choice of distribution $G$; when the latter is involutive, the description simplifies, see  Corollary \ref{rem:infdefcodim1}. 

\subsubsection*{Space-filling branes}

 For the reader's convenience, we first explain  the space-filling case.
 We denote by $\Omega^{(1,1)}(M,\CC)$ the forms of type $(1,1)$ w.r.t. the complex structure $I:=\omega^{-1}\circ F$; they are characterized by $B(Iv,Iw)=B(v,w)$ for all vectors in the complexified tangent bundle.
 We denote by $\Omega^{(1,1)}_{\RR}(M)$ the real $(1,1)$ forms on $M$, i.e.\ the elements of $\Omega^{(1,1)}(M,\CC)$ invariant under conjugation.

Let $F$ be a space-filling brane structure on $(M,\omega)$.

\begin{lem}\label{lem:infdefspace-filling}
Let $\{F_t\}\subset \Omega^2(M)$ be a smooth curve with $F_0=F$ and $(\omega^{-1}\circ F_t)^2=-Id$ for all $t$. Then
$$\dot{F}:=\frac{d}{dt}|_0 F_t\in \Omega^{(1,1)}_{\RR}(M).$$
Further, if each $F_t$ is a space-filling brane, then $\dot{F}$ is closed.
\end{lem}
\begin{proof}
We have
$$0=\frac{d}{dt}|_0(\omega^{-1}\circ F_t)^2=\omega^{-1}\circ(\dot{F}\circ I + I^*\circ \dot{F}).$$ Since $\omega^{-1}$ is invertible, this shows that $\dot{F}$ is of type $(1,1)$.
\end{proof}
 
 \begin{cor}\label{cor:11closed}
Let $F$ be a space-filling brane. The infinitesimal deformations of $F$ (i.e.\ the formal tangent space at $F$ to the space of space-filling branes) are given by the closed elements of $\Omega^{(1,1)}_{\RR}(M)$.  
 \end{cor}

\subsubsection*{Arbitrary branes}

Recall from Definition \ref{def:brane} that a brane is a coisotropic submanifold $Y$ of $(M,\omega_M)$ together with a presymplectic form $F$, such that $F$ and $\omega:=\iota^*_Y\omega_M$ have the same constant rank kernel $E$, and so that on $TY/E$ the endomorphism $I:=\omega^{-1}\circ F$ satisfies $I^2=-Id_{TY/E}$.

Choose a complement $G$ to $E$. Applying Gotay's theorem (see Remark \ref{rem:gotay}) to the coisotropic submanifold $Y$, we may assume that $M$ is a neighborhood of the zero section in the  total space of the  vector bundle $q\colon E^*\to Y$. 

Consider a curve of branes through $(Y,F)$, given by $(Y_t,F_{Y_t})$, where $Y_t$ is the image of a section of $q$ and $F_{Y_t}$ is a closed 2-form there. 
This curve is encoded by a pair $(s_t,F_t)$ where 
$$s_t\in \Gamma(E^*) \text{   and   }F_t\in \Omega^2(Y),$$
 via $Y_t:=graph(s_t)$ and $(s_t)^*F_{Y_t}=F_t$. (In particular, $s_0(Y)=Y$ and $F_0=F$.) The pullback of $\omega_M$ to $Y_t$, once we identify the latter with $Y$, is \begin{equation}\label{eq:OPomegat}
     \omega_t:=\omega-d\bar{s
}_t.
\end{equation}
{Here, for any section $\sigma \in \Gamma(E^*)$, we denote by   $\bar{\sigma}$  the unique $1$-form on $Y$ such that  $\bar{\sigma}|_E=\sigma$ and $\bar{\sigma}|_G=0$ \cite[Corollary 3.2]{OP}.}

	\begin{figure}[h!]
		\begin{center}
		\begin{tikzpicture}[scale=1]

            \draw[thick] (-3,0) -- (3,0); 
            		\node[right] at (3,0) {$Y$};
          \draw (-3,-1)-- (-3,2);
          \draw (-2,-1)-- (-2,2);
          \draw (-1,-1)-- (-1,2);
          \draw (0,-1)-- (0,2);
          \draw (1,-1)-- (1,2);
            \draw (2,-1)-- (2,2);
              \draw (3,-1)-- (3,2); 
              \node[right] at (3,2) {$E^*$};  
              
              \draw[thick,cyan] (-3,1.5) .. controls (-1,-0.5) and (1,2) .. (3,1);
  \node[left] at (-3,1.5)[cyan] {$Y_t=graph(s_t)$};  
			\end{tikzpicture}
		\end{center}
		\caption{The vector bundle $q\colon E^*\to Y$ with a section $s_t$.}
	\end{figure}

We now take the time derivative at $0$ of the brane condition, and obtain the   constraints below, using the notation
$$\dot{s}:=\frac{d}{dt}|_0 s_t,\quad \dot{F}:=\frac{d}{dt}|_0 F_t,\quad\dot{\omega}:=\frac{d}{dt}|_0\omega_t.$$
 We will be using the following terminology: a 2-form on the coisotropic submanifold $Y$ is \emph{horizontal} if it vanishes on  $\wedge^2 E$, i.e.\ if its pullback to the leaves of the characteristic distribution vanishes.

\begin{lem}\label{lem:infdefdirect}
Taking the time derivative at $0$ of the brane condition for the curve $\{(Y_t,F_{Y_t})\}$, one obtains exactly the following linear conditions on $\dot{s}$ and $\dot{F}$:

\begin{itemize}
    \item[i)]  The 2-form  $\dot{\omega} $  on $Y$ satisfies  $\dot{\omega}=-d{\bar{\dot{s}}}$ and is horizontal.    
        \item[ii)] The 2-form  $\dot{F}$ on $Y$ is  closed and horizontal.
    \item[iii)]   For any $e\in E$, the covectors $\dot{F}(e)$ and $\dot{\omega}(e)$ (which by the above items lie in the annihilator 
  $ E^{\circ} \cong (TY/E)^*$) satisfy:
    \begin{equation}\label{eq:samekernelDirac}       \dot{F}(e)=I^*(\dot{\omega}(e)).
    \end{equation}

  \item[iv)]  For all $X\in TY$ we have the following equality of elements of  $E^{\circ}$:
  \begin{equation}\label{eq:I^2Dirac}
  \dot{\omega} (X)+  \dot{F} (\widehat{I[X]}) =I^*(\dot{\omega}\big( \widehat{I[X]})-\dot{F}(X)\big).
\end{equation}
Here $[X]$ denotes the image of $X$ under $TY\to TY/E$, and  $\widehat{I[X]}$ denotes any\footnote{Note that on the l.h.s. and r.h.s. of eq. \eqref{eq:I^2Dirac}, the same preimage of $X$ occurs.}
 preimage of   $I[X]$.
\end{itemize}
\end{lem}

\begin{rmk}[Interpretation of Lemma \ref{lem:infdefdirect}]\label{rem:implied}
Notice that $\dot{F}$ and $\dot{\omega}=-d{\bar{\dot{s}}}$ lie in $\Gamma(\wedge^2G^*)\oplus \Gamma(E^*\otimes G^*)$, by i) and ii) in the above lemma.

Item iii) is a condition on the components in  $E^*\otimes G^*$ of these two forms. Assuming item iii), the next item iv) is a condition on their $\wedge^2G^*$ components, since when $X\in E$ it reduces to item iii). {(Even though condition iii) is implied by condition iv), we prefer to phrase  Lemma \ref{lem:infdefdirect} using both, as they are conditions on two different components).}
\end{rmk}

 \begin{rmk}[On time independence]\label{rem:omegatindep}
 
a) When $\omega_t$ is time-independent,  
 eq. \eqref{eq:samekernelDirac} in item iii) boils down to  $E\subset \ker \dot{F}$, which is a stronger condition than  $\dot{F}$ being horizontal (the condition in ii)).

b) Assume that the restriction\footnote{Here 
$\omega_t|_G$ denotes 
the restriction of $\omega_t$ to a bilinear form on $G$.} $\omega_t|_G$ is time-independent.
 Then
eq. \eqref{eq:I^2Dirac}  in item iv) boils down to the condition of item iii) and  to $\dot{F}|_G I+I^*\dot{F}|_G$=0, the latter being just the condition that $\dot{F}|_G$ is a form of type $(1,1)$.
Here 
we view $I$ as an endomorphism of $G$ via  $G\cong TY/E$.  

Notice that the assumption on $\omega_t|_G$ being time-independent holds whenever $G$ is involutive, since in that case 
$\bar{\dot{s}}|_G=0$ implies
$(d\bar{\dot{s}})|_G=0$, hence $\dot{\omega}|_G=0$ by item i).
\end{rmk}

In the following remark we outline how the conditions of Lemma \ref{lem:infdefdirect} correspond to the various requirements in the definition  of brane.
The proof of Lemma \ref{lem:infdefdirect}, which in particular shows that no further conditions arise, is deferred to   Appendix \ref{app:proof}.

\begin{rmk}
\emph{Item i)} follows from the coisotropic condition on $Y_t$, and \emph{item ii)}   from the fact that the $F_t$ are presymplectic forms on $Y$ (i.e.\ closed and of constant rank).

\emph{Item iii)} follows from the fact that    $F_t$ and $\omega_t$ have the same kernel, for each $t$. Indeed, {given a point $p\in Y$,}
 let $e_t$ be any smooth family   {in $T_pY$} with $e_t\in \ker(F_t)_p$ and $e_0=e$. 
 Then
\begin{equation}\label{eq:Fe}
0=\frac{d}{dt}|_0(F_t(e_t,\cdot))=\dot{F}(e,\cdot)+F(\dot{e},\cdot),
\end{equation} where $\dot{e}:=\frac{d}{dt}|_0e_t\in {T_pY}$.
By item ii),   $\dot{F}(e,\cdot)$ annihilates $E$, hence it induces an element of $(TY/E)^*$, denoted using square brackets. 
Equation \eqref{eq:Fe} explains the first equality in
$$[\dot{F}(e,\cdot)]=-[F(\dot{e},\cdot)]=-[\omega(\dot{e},\cdot)]\circ I=
[\dot{\omega}(e,\cdot)] \circ I,$$
whereas the middle equality follows from
 the definition of $I$.
For the third equality, notice that the analogue of \eqref{eq:Fe} holds for $\omega$ too, since  $F_t$ and $\omega_t$ have the same kernel.
 
\emph{Item  iv)}  follows from the fact that   $I_t:=(\omega_t)^{-1}\circ F_t$ squares to $-Id$.
Indeed, on any distribution   on $Y$ transverse to $E$, similarly to Lemma \ref{lem:infdefspace-filling}, one can compute
$\frac{d}{dt}|_0(I_t)^2$
  using $\frac{d}{dt}|_0(\omega_t)^{-1}=-\omega^{-1}\circ \dot{\omega}\circ\omega^{-1}$.
\end{rmk}

In condition i) of Lemma \ref{lem:infdefdirect}, the fact that $-d{\bar{\dot{s}}}$  is horizontal is equivalent to the foliated $1$-form $\dot{s}$ being closed. Hence Lemma \ref{lem:infdefdirect} yields the following,  by taking $B=\dot{F}$ and $r=\dot{s}$ there:

\begin{cor}\label{cor:infdefdirect}
Choose a  distribution  $G$ such that  $G\oplus E=TY$. 
The infinitesimal deformations of a brane $(Y,F)$
are given by pairs 
$$(r,B)\in \Gamma(E^*)\oplus \Omega^2(Y)$$
 such that  $r$ is  closed  as a foliated 1-form  and:
 \begin{itemize} 
        \item[ii)] The 2-form  $B$  is  closed and horizontal.
    \item[iii)]   For any $e\in E$, the covectors $B(e)$ and $  (d\bar{r})(e)$  in 
  $ E^{\circ}$ satisfy:
    \begin{equation}\label{eq:samekernelDirac2}
        B(e)=-I^*(d\bar{r}(e)),
    \end{equation}
 where    $\bar{r} \in \Omega^1(Y)$ denotes the unique extension of $r$ that annihilates $G$. 

  \item[iv)]  For all $X\in TY$ we have the following equality of elements of  $E^{\circ}$:
  \begin{equation}\label{eq:I^2Dirac2}
-d\bar{r} (X)+  B (\widehat{I[X]}) =I^*(-d\bar{r}\big( \widehat{I[X]})-B(X)\big).
\end{equation}
Here $[X]$ denotes the image of $X$ under $TY\to TY/E$, and  $\widehat{I[X]}$ denotes any
 preimage of   $I[X]$.
\end{itemize}
\end{cor}

{Notice that condition iii) is implied by condition iv), see Remark \ref{rem:implied}.
}

\begin{rmk}\label{rem:codim1}
i) It is known \cite{OP} that the infinitesimal deformations of $Y$ as a coisotropic submanifolds are given by $\Gamma_{\text{cl}}(E^*)$, the closed foliated 1-forms on $Y$. This fact is recovered by  Corollary \ref{cor:infdefdirect}.

ii) Notice that when the submanifold $Y$ has codimension $1$, the conditions that  $r$ is closed and that $B$ is horizontal are automatically satisfied. 
\end{rmk}

\begin{rmk}[Dependence on the  choice of complement $G$]
 The concrete description of infinitesimal deformations given in Corollary \ref{cor:infdefdirect} depends on the choice of complement $G$. For given $r\in \Gamma(E^*)$, the extension $\bar{r} \in \Omega^1(Y)$  depends on $G$. 
 Further, the r.h.s of \eqref{eq:samekernelDirac2} really depends on the choice of $G$. Indeed,
 given tangent vectors $e\in E_p$ and $v\in G_p$ at $p\in Y$, one computes that $$(\iota_e d\bar{r})(v)=-(\cL_{\tilde{v}}r)(e)$$ where $\tilde{v}\in \Gamma(G)$ is an   extension of $v$ in a neighborhood $U\subset Y$ to a vector field which is projectable w.r.t. the quotient map $U\to U/E$.
 If $G_1$ and $G_2$ are two complements of $E$, $\bar{r}_i$ the corresponding extensions, $v_i\in (G_i)_p$ defining the same element in $T_pY/E_p$, and the extensions $\tilde{v_i}$ are chosen to project to the same vector field on $U/E$,
 then the difference $(\iota_e d\bar{r}_1)(v_1)- (\iota_e d\bar{r}_2)(v_2)$ equals $(d_E(\iota_{\tilde{v}_1-\tilde{v}_2}r))(e)$, which is generally non-zero.  {Here $d_E$ is the foliated differential.} 
\end{rmk}

\subsubsection*{Transversely integrable branes}
 
 Condition iv) in Corollary \ref{cor:infdefdirect} simplifies when the complement distribution $G$ can be chosen to be involutive, by   Remark \ref{rem:omegatindep} b), yielding the following.

\begin{cor}
\label{rem:infdefcodim1} Let $(Y,F)$ be a  brane. Suppose there is an  \emph{involutive} distribution  $G$ such that  $G\oplus E=TY$.
Then    the  {infinitesimal}  deformations of $Y$ as a brane  are   given by pairs  $$(r,B)\in \Gamma(E^*)\oplus \Omega^2(Y)$$ such that  $r$ is closed as a foliated 1-form and 
\begin{itemize}
\item[ii)] The 2-form  $B$  is  closed and horizontal.
\item[iii)] for all $e\in E$, the covectors $B(e)$ and $ (d\bar{r})(e)$  in   $ E^{\circ}$ satisfy:
\begin{equation}\label{eq:mixedcondition}
  B(e)=-I^*(d \bar{r}(e)),
\end{equation}
where $\bar{r}\in \Omega^1(Y)$ is the extension of $r$ annihilating $G$
\item[iv)] $B|_G$ is of type $(1,1)$, i.e.\ $B(IX,IX')=B(X,X')$ for all $X,X'\in G$.
\end{itemize}
\end{cor}

 Notice that ii), iii) and iv) are conditions on the restriction of $B$ to $\wedge^2 E$, to $E\otimes G$ and to $\wedge^2 G$ respectively.

\begin{ex}[Lagrangian branes]
 Let $Y$ be a Lagrangian submanifold of $(M,\omega_M)$. We saw in Example \ref{ex:lagr} that the zero 2-form  makes $Y$ into a brane, and its kernel is $E=TY$.
  Thus in Corollary \ref{rem:infdefcodim1} we have $B=0$, and we recover the fact that infinitesimal deformations of Lagrangian submanifolds are given by closed $1$-forms on $Y$.  
\end{ex}

\begin{ex}[Space-filling branes]\label{rmk:defsspacefill}
For space-filling branes  $Y=M$,
 in Corollary \ref{rem:infdefcodim1} we have 
$E=0$. Hence  infinitesimal deformations  are given by
closed elements of $\Omega^{(1,1)}_{\RR}(M)$, recovering Lemma \ref{lem:infdefspace-filling}. 
\end{ex}

\begin{ex}[Space-filling branes on $\RR^4$]\label{ex:4torusbranes}
Consider  $M=\CC^2\cong \RR^4$, let $F$ and $\omega$ be the real and imaginary part of $dz_1\wedge dz_2$, as in Example \ref{rem:IdfcoordsA}.

We can write out a frame for $\Omega^{(1,1)}_{\RR}(M)$ by taking the real and imaginary parts of $dz_1\wedge d\bar{z}_1$, $dz_2\wedge d\bar{z}_2$ and $dz_1\wedge d\bar{z}_2$, which are
$$dx_1\wedge dy_1,\quad   dx_2\wedge dy_2, \quad  dx_1\wedge dx_2+ dy_1\wedge dy_2, \quad 
 -dx_1\wedge dy_2+dy_1\wedge dx_2.$$
As seen in Example \ref{rmk:defsspacefill}, the infinitesimal deformations of the space-filling brane structure $F$ on $(M,\omega)$ are given by the \emph{closed} 
$C^{\infty}(M)$-linear combinations of these four differential forms. The same holds for the 4-torus $\RR^4/\ZZ^4$.
\end{ex}

For the codimension 1 branes introduced in Example \ref{ex:mainex}, the infinitesimal deformations will be described explicitly in 
Lemma \ref{lem:infdefmainex}.
 
 \begin{rmk}[Comparison with the literature]\label{complit}
We compare briefly the infinitesimal deformations we obtained in this subsection with those that appeared in the work of Koerber-Martucci \cite{KoerMar} (and later in Collier's work \cite{Collier}). For branes $(Y,F)$ in arbitrary generalized complex manifolds in the standard  Courant algebroid  (so  $H\in \Omega^3(M)$ vanishes), the infinitesimal deformations  they consider are of the form  $(r,B)\in \Gamma(E^*)\oplus \Omega^2(Y)$ with $B$ an \emph{exact} 2-form \cite[\S 3, page 9]{KoerMar}.

As the 2-form  $B$ in Corollary \ref{cor:infdefdirect} is closed but not necessarily exact, the infinitesimal deformations we consider are more general. For instance, on the 4-torus as in Example
\ref{ex:4torusbranes}, the infinitesimal deformation $B:=dx_1\wedge dy_1$ is not exact. Notice that there is a 1-parameter family of space-filling branes
$F_t:=F+tB$, where $t\in \RR$
(see \cite[Lemma 4.3 and 4.4]{KLZTorelli} or  check directly that $(\omega^{-1}\circ F_t)^2=-Id$), for which the two-form $B$ is the initial velocity.

From the description for branes in symplectic manifolds given in \cite[\S 7.3.2]{KoerMar}, it appears that the infinitesimal deformations that are regarded as  trivial there (see equation (7.45) there) agree with those that we will regard as trivial in \S \ref{subsec:cochaincomplexes} (see in particular Remark \ref{rem:spacefillfunctions} for space-filling branes).
 \end{rmk}

\subsection{Cochain complexes governing infinitesimal deformations of   branes}\label{subsec:cochaincomplexes}

In \S \ref{subsec:Diracapproach} we described the infinitesimal deformations of a brane $(Y,F)$. Now we display cochain complexes which \emph{control}  the infinitesimal deformations. By this we mean that
\begin{itemize}
    \item 
the 1-cocycles correspond to infinitesimal deformations of $(Y,F)$, and 
\item the 1-coboundaries correspond to infinitesimal deformations that arise from hamiltonian isotopies (and therefore might be regarded as trivial). 
\end{itemize}
In particular, denoting by $H^1$ the first cohomology of the complex, we will have
$$H^1\cong \frac{ \text{Infinitesimal deformations}}{\text{Infinitesimal deformations by hamiltonian  symmetries}}=T_{(Y,F)}\cM,$$
where the latter is  the formal tangent  space to the space of  branes modulo hamiltonian isotopies.

As in \S \ref{subsec:Diracapproach}, we start considering space-filling branes, then we pass to arbitrary branes, and finally we spell out the special case in which the characteristic distribution admits an involutive complement.

\subsubsection*{Space-filling branes}
Let {$F$ be a space-filling brane structure} on the symplectic manifold $(M,\omega)$.

\begin{rmk}[Hamiltonian deformations]\label{rem:spacefillfunctions}
    Let $f\in C^{\infty}(M)$. The  flow $\Phi_t$ of the  hamiltonian vector field $X_f$  preserves $\omega$. Hence
    applying the flow to $F$    we obtain a family of space-filling brane $(M, \Phi_{-t}^*F)$ in  $(M,\omega)$. Notice that
    \begin{equation}\label{eq:did}
        \frac{d}{dt}|_0 \Phi_{-t}^*F=-\cL_{X_f}F=-d\iota_{X_f}F=-dI^*df,
    \end{equation}
    using $F=I^*\circ \omega$ in the last equality.
    Hence for any smooth function $f$, we see that $-\cL_{X_f}F$ is an infinitesimal deformation of the space-filling brane $(M,F)$. Notice that 
$-\cL_{X_f}F$ is a $(1,1)$-form: this is a consequence of Lemma \ref{lem:infdefspace-filling}, but can be also seen directly from $dI^*df=-2i \del\delb f$.
\end{rmk}

 Lemma \ref{lem:infdefspace-filling} and Remark \ref{rem:spacefillfunctions} imply the following statement.
\begin{prop}\label{prop:cplxspace-filling}
    Let $(M,F)$ be a space-filling brane. Then the following cochain complex controls the infinitesimal deformations of the brane:
    \begin{equation}\label{eq:complexspace}
        C^{\infty}(M)\longrightarrow  \Omega^{(1,1)}_{\RR}(M) \longrightarrow
    \Omega^3(M) \longrightarrow \dots
    \end{equation}
    where the first map is $f\mapsto \cL_{X_f}F$ and the other maps are given by the de Rham differential.
\end{prop}

\begin{rmk}[The first cohomology]
The first cohomology of the complex \eqref{eq:complexspace} depends only on the complex structure $I$, as is clear from \eqref{eq:did}. Even more,  it agrees with the real part of the first Bott-Chern cohomology group 
of the complex manifold $(M,I)$, by the very definition of the latter.
\end{rmk}

\begin{rmk}\label{rem:cpctK}[The compact Kähler case]
Suppose that $(M,I)$ is compact Kähler, where $I$ is the complex structure determined as $I:=\omega^{-1}\circ F$.
 The infinitesimal deformations of $F$ by hamiltonian vector fields considered in Remark \ref{rem:spacefillfunctions} are precisely the $d$-exact elements of $\Omega^{1,1}_{\RR}(M)$, as a consequence of Lemma \ref{rem:additionalnote} below.
 Therefore, in the compact Kähler case, 
 we have $T_{(Y,F)}\cM\cong  H^{(1,1)}_{\RR}(M)$.
The latter is defined as the closed real $(1,1)$ forms modulo the exact ones, and is isomorphic to the real part of the Dolbeaut cohomology in degree $(1,1)$, see  \cite[Lemma A.2]{KLZTorelli}. 
This remark can also be deduced from the general fact that, on a compact Kähler manifold,  
Bott-Chern cohomology is isomorphic to Dolbeaut cohomology.
\end{rmk}

\begin{lem}\label{rem:additionalnote}
On a compact Kähler manifold $(M,I)$, the
$d$-exact real $(1,1)$-form are given  precisely by $\{d(I^* d f): f\in C^{\infty}(M)\}$.   

Further, the latter equals $\{d(I^* \alpha): \alpha\in \Omega^1_{\text{cl}}(M)\}$.
\end{lem}
\begin{proof}
     Recall that $d^c$ is the real operator on complex differential forms given by $d^c=-i(\del-\delb)$. It satisfies $dd^c=-d^cd=2i\del\delb$.
 The $d^cd$-Lemma \cite[Remark 3.A.23]{Huy}  states the following on a compact K\"ahler manifold: if $\alpha$ is a $d^c$-closed and $d$-exact complex $k$ form, then $\alpha=d^cd\beta$ for some  complex $k-2$ form.
 
 Now let $\alpha\in \Omega^{1,1}_{\RR}(M)$ be $d$-exact. Then $\alpha$ is $d$-closed, and comparing bi-degrees we see that $\del \alpha =0=\delb \alpha$. Thus  $\alpha$ is also $d^c$-closed. Hence by the  $d^cd$-Lemma there is a function $f$ (which we may take to be real valued) such that $\alpha=d^cdf=-dd^cf=d I^*df$. Here we used $I^*df=i\del f - i \delb f=-d^cf$.

 The converse implication is shown by a direct computation, valid on any complex manifold: for any closed $1$-form $\alpha$, one has $d(I^*\alpha)\in \Omega^{(1,1)}_{\RR}(M)$. This also proves the last statement.
\end{proof}

   In Proposition \ref{prop:cplxspace-filling}   we saw that
\eqref{eq:complexspace} is a cochain complex controlling the infinitesimal deformations of a space-filling brane.
There are several complexes
controlling such deformations, and 
we now argue that \eqref{eq:complexspace} is a good choice, because of the 
 following  result on obstructions.
Denote by  $F^{-1}$ the Poisson bivector field  inverse to the symplectic form $F$, and by $[\;,\;]_{F^{-1}}$ the associated Koszul bracket: 
 on 1-forms it is defined by
$$[\alpha,\beta]_{F^{-1}}=\cL_{F^{-1}\alpha}\beta-\cL_{F^{-1}\beta}\alpha - d(F^{-1}(\alpha,\beta)),$$
and it is extended
to arbitrary differential forms as a graded biderivation of the wedge product
(see e.g.\ \cite{PoisGeoBookAMS}).
We claim that if an infinitesimal brane deformation $B$ is not mapped to zero under the map
\begin{align}\label{eq:kura}
\Omega^{(1,1)}_{\RR,\text{cl}}(M)&\to 
\Omega_{\text{cl}}^{3}(M)/ d \Omega^{(1,1)}_{\RR}(M),\\
B&\mapsto 
[B,B]_{F^{-1}} \;\text{ mod }d \Omega^{(1,1)}_{\RR}(M),
\nonumber
\end{align}
then $B$ can not be prolonged, i.e.\ there exists no smooth 1-parameter family $F_t$ of space-filling branes such that $F_0=F$ and $\frac{d}{dt}|_0 F_t=B$. 

 In the following remark we justify the above  obstructedness result;
details and proofs will appear in a separate paper.

\begin{rmk}[On the obstruction map \eqref{eq:kura}]
Let $(M,\omega)$ be a symplectic manifold, and $F$ a space-filling brane structure. In particular, $F$ is a symplectic form, hence\footnote{To take the graph we view $F$ as a map $TM\to T^*M, v\mapsto \iota_vF$.}   $\graph(F)$ is a Dirac structure \cite{Cou}. To describe the deformations of $\graph(F)$ as a Dirac structure, one can use an auxiliary complementary Dirac structure:
the latter provides a parametrization of  Dirac structures nearby $\graph(F)$ by elements of $\Gamma(\wedge^2 (\graph(F))^*)$ satisfying the Maurer-Cartan equation
of an induced differential graded Lie algebra (DGLA) structure on $\Gamma(\wedge (\graph(F))^*)[1]$ \cite{LWX}.
Under the identification $\graph(F)\cong TM$ by the first projection, the above space of sections is just $\Omega(M)[1]$, where ``$[1]$'' denotes a degree shift.

It turns out that a good choice of complement is $\graph(-F)$; the DGLA structure inherited by $\Omega(M)[1]$ is then given by the de Rham differential $d$ and the   Koszul bracket of $-\frac{1}{2}F^{-1}$.
The reason why $\graph(-F)$ is a good choice of complement is that  under the associated parametrization of Dirac structures, the graph of a closed 2-form is a space-filling brane structure if and only if it lies in a \emph{linear} subspace   of $\Omega^2(M)$, namely $\Omega^{(1,1)}_{\RR}(M)$. It follows that the sub-DGLA $$(\Omega^{(1,1)}_{\RR}(M)\oplus_{k \ge 3}\Omega^{k}(M)[1], d, [\;,\;]_{-\frac{1}{2}F^{-1}})$$
has the property that its Maurer-Cartan elements parametrize
 the deformations of the space-filling brane structure $F$.
The general theory of deformations via DGLA's (see e.g.\ \cite[Theorem 11.4]{OP}) then implies the obstructedness result just before this remark.

 We are currently investigating how include in this picture the Hamiltonian deformations mentioned in Remark \ref{rem:spacefillfunctions}.
\end{rmk}

\begin{ex}[Infinitesimal deformations of space-filling branes on $\TT^4$] \label{ex:TT4obstrmap}
Consider the $4$-torus $\TT^4$, endowed\footnote{A concrete instance of this setting was given in Example \ref{ex:4torusbranes}.} with a symplectic form $\omega$ admitting a space-filling brane $F$. Then the   obstruction map \eqref{eq:kura} vanishes identically.

Indeed, we may assume that  $\omega$, $F$ and $I$ are all constant on the torus, since the first claim in the proof of \cite[Prop. 4.17]{KLZTorelli} states that we can always find a diffeomorphism of the torus yielding this. Since the Poisson bivector field  $F^{-1}$ is constant,  the associated Koszul bracket vanishes on constant elements of $\Omega^{(1,1)}_{\RR,\text{cl}}(\TT^4)$. 

An arbitrary element  $B\in \Omega^{(1,1)}_{\RR,\text{cl}}(\TT^4)$ is of the form $B_0+d\alpha$ where $B_0$ is constant $(1,1)$ form and $\alpha$ is a 1-form. The Koszul bracket 
$[B_0+d\alpha,B_0+d\alpha]_{F^{-1}}$ is exact, as it equals $2d[\alpha,B_0]_{F^{-1}}+d[\alpha,d\alpha]_{F^{-1}}$. Notice that $\Omega^{3}(\TT^4)=\Omega^{(2,1)+(1,2)}_{\RR}(\TT^4)$ because of dimension reasons.
Now, every exact element of $\Omega^{(2,1)+(1,2)}_{\RR}(\TT^4)$ lies in $d \Omega^{(1,1)}_{\RR}(\TT^4)$, implying that the image of $B_0+d\alpha$ under the obstruction map \eqref{eq:kura} vanishes.

The last statement holds because on every compact K\"ahler manifold $M$, any exact element $C\in \Omega^{(2,1)}(M,\CC)\oplus \Omega^{(1,2)}(M,\CC)$ lies in $d\Omega^{(1,1)}(M,\CC)$. Indeed, take any primitive $\epsilon=\epsilon^{(2,0)}+\epsilon^{(1,1)}+\epsilon^{(0,2)}$ of $C$. Since $\bar{\partial}\epsilon^{(2,0)}$ is $d$-closed, by the the $\partial \bar{\partial}$-Lemma \cite[Remark 3.2.10]{Huy} we have $\bar{\partial}\epsilon^{(2,0)}=\partial \gamma$ for some $\bar{\partial}$-closed $(1,1)$-form $\gamma$. Similarly, ${\partial}\epsilon^{(0,2)}=\bar{\partial} \delta$ for some $\partial$-closed $(1,1)$-form $\delta$. Then $\epsilon^{(1,1)}+\gamma+\delta \in \Omega^{(1,1)}(M,\CC)$ 
is also a primitive of $C$.
\end{ex}

Since the map \eqref{eq:kura} provides an obstruction to the prolongability to infinitesimal deformations, Example \ref{ex:TT4obstrmap} is consistent with the next lemma, whose proof we defer to the appendix. 

\begin{lem}\label{lem:T4unobstr}
Consider the $4$-torus $\TT^4$, endowed with a symplectic form $\omega$ admitting a space-filling brane $F$. Then any 
infinitesimal deformation $B$ of $F$ can be prolonged to a smooth 1-parameter family of space-filling branes.
\end{lem}

\subsubsection*{Arbitrary branes}

Now   let $(Y,F)$ be any  brane in $(M,\omega_M)$.
Choose a complement $G$ to the characteristic distribution $E$. Applying Gotay's theorem (see Remark \ref{rem:gotay}) to the coisotropic submanifold $Y$, we may assume that $M$ is a neighborhood of the zero section in the   total space of the  vector bundle $q\colon E^*\to Y$.
We denote by   $\omega$ the pullback of $\omega_M$ to $Y$.

\begin{rmk}[Hamiltonian deformations]\label{rmk:branefct}
        Let $f\in C^{\infty}(Y)$. Consider the  flow $\Phi_t$ of the  hamiltonian vector field $X_{q^*f}$ w.r.t. $\omega_M$. 
    Applying the flow to $(Y,F)$, {for values of $t$ sufficiently close to zero}    we obtain a family of  branes supported on $\graph(s_t)$ where $s_t\in \Gamma(E^*)$.
The corresponding $2$-form, under the identification $\graph(s_t)\cong Y$ provided by the projection $q$, is
$$F_t:=s_t^*(\Phi_t^{-1})^*F=[(q\circ \Phi_t)^{-1}]^*F.$$
Notice that {when $t$ is sufficently small,}  $(q\circ \Phi_t)^{-1}$ is a 1-parameter family of diffeomorphisms of $Y$ satisfying
    \begin{equation}
        \frac{d}{dt}|_0 (q\circ \Phi_t)^{-1}|_Y=
        -\frac{d}{dt}|_0 (q\circ  \Phi_t)|_Y=
        -q_*(X_{q^*f}|_Y)=-X_f
    \end{equation}
where  $X_f:=(\omega|_G)^{-1}(df|_G)$.

In the last equation we used that $\omega_M$ at points of $Y$  is given by the sum of $\omega|_G$ and the canonical skew pairing between $E^*$ and $E$ \cite[Eq. (4)]{Gotay}, and consequently\footnote{{Here $d_E$ is the foliated differential.}}  $X_{q^*f}|_Y=X_f-d_Ef\in \Gamma( TY\oplus E^*)$.
Hence
$$\frac{d}{dt}|_0 F_t=-\cL_{X_f}F.$$
Further, $\frac{d}{dt}|_0 s_t=-d_Ef\in \Gamma(E^*)$, being this the vertical part of $X_{q^*f}|_Y$.

We conclude that the pair  $(-d_Ef,-\cL_{X_f}F)$ is an infinitesimal deformation of the space-filling brane $(M,F)$.  
\end{rmk}

For all integers $k$ define
$$\Omega^{k}_{\hor}(Y):=\left\{ B\in \Omega^k(Y) \text{ s.t.  $B|_{\wedge^k E}=0$}\right\},$$
the subcomplex of the de Rham complex given by horizontal forms.

\begin{prop}\label{prop:defcomplxinv2}
    Let $(Y,F)$ be a  brane. Fix  a distribution  $G$ such that  $G\oplus E=TY$. Then the following cochain complex controls the infinitesimal deformations of the brane: 
    \begin{equation}\label{eq:defcplxinv2}
        C^{\infty}(Y)\longrightarrow 
        \left\{(r,B)\in\Gamma(E^*)\oplus \Omega^{2}_{\hor}(Y):  
 \text{eq. \eqref{eq:samekernelDirac2} and \eqref{eq:I^2Dirac2}     hold}   
        \right\}
\longrightarrow
        \Gamma(\wedge^2 E^*)\oplus \Omega^3_{\hor}(Y) \longrightarrow \dots
    \end{equation}
    where the first map is $f\mapsto (d_Ef,\cL_{X_f}F)$ for    $X_f:=(\omega|_G)^{-1}(df|_G)$, and the other maps are $d_E\oplus d$.      
\end{prop}
{Recall that eq. \eqref{eq:samekernelDirac2} is redundant as it is  implied by  eq. \eqref{eq:I^2Dirac2}, see the text under Corollary \ref{cor:infdefdirect}.}

\begin{proof}
The first map  is well-defined, by  Corollary \ref{cor:infdefdirect}  and Remark \ref{rmk:branefct}. The above is a cochain complex (use that $F$ is closed).  The infinitesimal deformations of the brane are precisely the closed elements of the second space appearing in \eqref{eq:defcplxinv2}, as follows again using Corollary \ref{cor:infdefdirect}.
\end{proof}

\begin{lem}\label{lem:equivprolong}
Let $(Y,F)$ be a  brane and $G$  a distribution complementary to $E$.
Fix an infinitesimal deformation $(r,B)$, i.e.\ a degree $1$ cocycle in the cochain complex \eqref{eq:defcplxinv2}.
Suppose that $(r,B)$  can be prolonged to a smooth curve of branes.
 Then any infinitesimal deformation lying in the same cohomology class of $(r,B)$ can also be prolonged to a smooth curve of branes.
\end{lem}
\begin{proof}
By assumption, there is a smooth family
 $(s_t,F_t)\in \Gamma(E^*)\oplus \Omega^2(Y)$ so that $\frac{d}{dt}|_0(s_t,F_t)=(r,B)$
 and so that $(\graph(s_t), (q|_{graph(s_t)})^*F_t)$ is a family of branes, see \S \ref{subsec:Diracapproach}.

Any infinitesimal deformation lying in the same cohomology class of $(r,B)$ can be written as $(r-d_Ef,B-\cL_{X_f}F)$ for some $f\in C^\infty(Y)$.

Denote by $\Psi_t$ the flow of the hamiltonian vector field $X_{q^*f}$.
Since $\Psi_t$ is a symplectomorphism, we see that $(\Psi_t(\graph(s_t)), (\Psi_t^{-1})^*(q|_{graph(s_t)})^*F_t)$ is a family of branes.
For $t$ sufficiently small, 
we can write it as
$$(\graph(\sigma_t), (q|_{\graph(\sigma_t)})^*\Phi_t)$$
for a smooth family $(\sigma_t,\Phi_t)\in \Gamma(E^*)\oplus \Omega^2(Y)$.

We claim that 
\begin{equation}\label{eq:claim1}
    \frac{d}{dt}|_0 \sigma_t=\frac{d}{dt}|_0 s_t-d_Ef.
\end{equation}
Indeed, denote   $\psi_t:=q\circ \Psi_t\circ s_t$, and let $y\in Y$. Notice that
$\frac{d}{dt}|_0 (\Psi_t\circ s_t)(y)=X_{q^*f}(y)+\frac{d}{dt}|_0s_t(y),$ so 
\begin{equation}\label{eq:derpsi}
    \frac{d}{dt}|_0 \psi_t(y)=X_{q^*f}(y)_{hor},
    \end{equation}
where the index denotes the first component in the decomposition $T_yE^*=T_yY\oplus E^*_y$. By construction, $(\sigma_t\circ \psi_t)(y)=(\Psi_t\circ s_t)(y)$. Taking the time derivative of this equation at $t=0$  we see that
$\frac{d}{dt}|_0 \sigma_t(y)=X_{q^*f}(y)- X_{q^*f}(y)_{hor}+\frac{d}{dt}|_0s_t(y)$.
Using that the vertical component of 
$X_{q^*f}(y)$ is $-d_Ef(y)$
(see Remark \ref{rmk:branefct}), this proves the claim.

Further, 
\begin{equation}\label{eq:claim2}
\frac{d}{dt}|_0\Phi_t=\frac{d}{dt}|_0F_t-\cL_{X_f}F,
\end{equation}
as one sees using $\Phi_t=(\psi_t^{-1})^*F_t$ and \eqref{eq:derpsi}.
Combining \eqref{eq:claim1} and \eqref{eq:claim2} finishes the proof.
\end{proof}

\subsubsection*{Transversely integrable branes}

For branes whose characteristic distribution admits an involutive complement $G$, the cochain complex of Proposition \ref{prop:defcomplxinv2} simplifies slightly, thanks to  Corollary \ref{rem:infdefcodim1}. Denote 
$$\Omega^{(1,1)}_{\hor}(Y):=
\left\{ B\in \Omega^{2}_{\hor}(Y) \text{ s.t.  $B|_{G}$ is of type $(1,1)$}\right\}.$$

\begin{prop}\label{prop:defcomplxinv}
    Let $(Y,F)$ be a  brane. Suppose there is an   \emph{involutive} distribution  $G$ such that  $G\oplus E=TY$. Then the following cochain complex controls the infinitesimal deformations of the brane: 
    \begin{equation}\label{eq:defcplxinv}
        C^{\infty}(Y)\longrightarrow 
        \left\{(r,B)\in\Gamma(E^*)\oplus \Omega^{(1,1)}_{\hor}(Y):  
     \text{eq. }\eqref{eq:mixedcondition} \text{ holds}  
        \right\}
\longrightarrow
        \Gamma(\wedge^2 E^*)\oplus \Omega^3_{\hor}(Y) \longrightarrow \dots
    \end{equation}
    where the first map is $f\mapsto (d_Ef,\cL_{X_f}F)$ for\footnote{So  $X_f$  is  the leafwise hamiltonian vector field of $f$ with respect to the foliation integrating $G$, endowed with the pullback of $\omega$ as symplectic form.}     $X_f:=(\omega|_G)^{-1}(df|_G)$, and the other maps are $d_E\oplus d$.   
\end{prop}

\begin{rmk}
 Assuming that $G$ is involutive, 
 we now prove in a direct way -- without using Remark \ref{rmk:branefct} nor  Corollary \ref{rem:infdefcodim1} -- that the first map in the complex \eqref{eq:defcplxinv} is well-defined. 
 Given a function $f$ on $Y$, denote 
 $$r:=d_Ef \quad \text{ and }B:=\cL_{X_f}F=d\iota_{X_f}F.$$
The  restriction of $B$ to the leaves of $G$ is of type $(1,1)$, as can be seen applying leafwise Remark \ref{rem:spacefillfunctions}. 

Let $e$ be a locally defined section of $E$, with the property that $\cL_e$ preserves $\Gamma(G)$
(i.e.\ so that $e$ is projectable under the quotient map to the local leaf space of the foliation integrating $G$). Then
$\iota_e B=\cL_e(\iota_{X_f}F)$ by Cartan's magic formula. Notice that this 1-form annihilates $E$ (since $E$ is involutive and is the kernel of $F$).
 This shows that $B$ is horizontal.
We now check \eqref{eq:mixedcondition}.
On $G$, we have
$$(\iota_e B)|_G=\cL_e(I^*(d_Gf))=I^* d_G(\cL_e f),$$
where the first equality holds as in Remark \ref{rem:spacefillfunctions}, and the second because $\cL_e$ preserves $I$ (since it preserves both of the closed forms 
$F$ and $\omega$). 
On the other hand, the 1-form $\bar{r}$ equals $d_Ef$ on $E$ and annihilates $G$, by definition. Hence for all $w\in \Gamma(G)$ we obtain $d\bar{r}(e,w)=-\cL_w(\bar{r}(e))=-\cL_w(\cL_e(f))$, making use of $[e,w]\in \Gamma(G)$. This shows that $\iota_e d\bar{r}=-d_G(\cL_e(f))$ (an equality of elements of $G^*\cong E^{\circ}$).
\end{rmk}

\begin{ex}[Lagrangian branes]
  For Lagrangian branes $Y$ we have $E=TY$, hence $\Omega^{\bullet}_{\hor}(Y)=\{0\}$.
The complex in Proposition \ref{prop:defcomplxinv} is therefore just the de Rham complex $(\Omega^{\bullet}(Y),d)$ of $Y$.
\end{ex}

\begin{ex}[Space-filling branes]\label{rmk:complexsspacefill}
For space-filling branes  we have $Y=M$  and hence 
$E=0$, in particular $\Omega^{\bullet}_{\hor}(Y)=\Omega^{\bullet}(Y)$. Hence the cochain complex \eqref{eq:defcplxinv} in Proposition \ref{prop:defcomplxinv} reduces to the cochain complex of Proposition \ref{prop:cplxspace-filling}.
\end{ex}

\subsection{Relation to infinitesimal deformations of coisotropic submanifolds}
\label{sec:infdefcoiso}

 Let $(Y,F)$ be a  brane in a symplectic manifold. We now consider the infinitesimal version of the forgetful map $\upsilon$  from branes to coisotropics displayed  in \eqref{eq:forgetful}, i.e.\ the formal derivative $\Upsilon$ of $\upsilon$ at $(Y,F)$. The non-surjectivity of $\Upsilon$ does not imply that  
$\upsilon$ is not locally surjective, i.e.\ it does not imply Theorem \ref{thm:nobrane}. However  when the {infinitesimal} deformations of the coisotropic are unobstructed (as in the codimension 1 case), the non-surjectivity of $\Upsilon$ implies that there is a \emph{smooth curve} of coisotropics that can not be lifted to a smooth curves of branes.

Recall that, upon choosing a  complement $G$ to the characteristic distribution $E$,  the vector space of  infinitesimal deformations of the brane  $(Y,F)$ is the vector subspace of $\Gamma(E^*)\oplus \Omega^2(Y)$ described in Corollary \ref{cor:infdefdirect}, and that the infinitesimal deformations of $Y$ as a coisotropic submanifold  are given by $\Gamma_{\text{cl}}(E^*)$.
{Notice that the first projection is a cochain map from the cochain complex \eqref{eq:defcplxinv2}  governing infinitesimal deformations of branes   to the cochain complex $(\Gamma(\wedge E^*),d_E)$ governing infinitesimal deformations of coisotropic submanifolds  \cite{OP};  in particular, it maps infinitesimal deformations to infinitesimal deformations.}
The main statement of this subsection is the following, which will be proven below by providing a counterexample.

\begin{prop}\label{prop:infnotsurj}
Let $(Y,F)$ be a  brane. The map 
\begin{equation}\label{eq:Upsilon}
  \Upsilon\colon 
\{\text{infinitesimal brane deformations}\}  
  \to \Gamma_{\text{cl}}(E^*),\;\;\;\;(r,B)\mapsto r
\end{equation}
is not surjective in general. 
\end{prop}

In other words, in general there are infinitesimal deformations of $Y$ as a coisotropic submanifold which do not arise from infinitesimal deformations
of $Y$ as a brane.
To prove this proposition we provide an explicit counter-example, see Proposition \ref{prop:notsurj} and Example \ref{rem:notsurj} below.

\begin{rmk}[The kernel of $\Upsilon$]\label{rem:kernelUps}
We have $$\ker(\Upsilon)\cong \{B\in \Omega^2(Y): B \text{ closed with $E\subset \ker(B)$ and $B|_G$ is of type $(1,1)$}\}.$$
Indeed, whenever $r=0$ one has $d\bar{r}=0$,  so the statement follows from Corollary \ref{cor:infdefdirect}, by similar arguments as in 
Remark \ref{rem:omegatindep}.
 When $\underline{Y}:=Y/E$ is a smooth manifold so that the projection $\pi \colon Y\to \underline{Y}$ is a submersion, we can describe $\ker(\Upsilon)$ simply as $\pi^*( \Omega^{1,1}_{\text{cl}}(\underline{Y}))$. Notice that the latter is isomorphic to the space of infinitesimal deformations of the space-filling brane structure $\underline{F}$ in the quotient $(\underline{Y},\underline{\omega})$, by
Corollary \ref{cor:11closed}. 
\end{rmk}

\subsubsection*{A codimension one example to prove  Proposition 
\ref{prop:infnotsurj}}
We now consider a particular class of branes, as in Example \ref{ex:mainex}. We use them to prove Proposition 
\ref{prop:infnotsurj}, via Proposition \ref{prop:notsurj} and Example \ref{rem:notsurj}.

\begin{prop}\label{prop:notsurj}
Fix a  symplectic manifold $(N,\omega_N)$ with a space-filling brane structure $F_N$.
We know from Example \ref{ex:mainex} that 
  $Y:=N\times S^1\times \{0\}$, together with the pullback $F$ of $F_N$, 
   is a brane
in the product symplectic manifold $(M:=N\times S^1\times\RR,\omega_M:=\omega_N\times \omega_{T^*S^1})$, with characteristic distribution $E=\{0\}\times TS^1$.

 The image of the map $\Upsilon$ in \eqref{eq:Upsilon} consists precisely of the elements $r\in\Gamma(E^*)$   such that
\begin{equation}\label{eq:imager}
   d_N I^* d_N \left(\int_{S^1}r\right)=0. 
\end{equation}
Here $I:=\omega_N^{-1}\circ F_N$, and
$\int_{S^1}r\in C^{\infty}(N)$ denotes the function obtained integrating $r$ along the fibers of $Y\to N$.
\end{prop}

\begin{rmk}
a) The domain of the map $\Upsilon$ depends on a choice of complement to $E=\{0\}\times TS^1$ of $Y$; here we choose the preferred (involutive) complement $G=TN\times \{0\}$.

b) The kernel of  the map $\Upsilon$ is the pullback to $Y$ of $\Omega^{1,1}_{\text{cl}}(N)$, see Remark \ref{rem:kernelUps}.
\end{rmk}


 To prove Proposition \ref{prop:notsurj}, we first give an explicit description of the infinitesimal deformations of the brane, in the next Lemma.

\begin{lem}\label{lem:infdefmainex}
 Assume the setting of Proposition \ref{prop:notsurj}.
Denote by $q$ the ``coordinate'' on $S^1=[0,1]/(0\sim 1)$, and write $N_q$ for $N\times \{q\}$.
The infinitesimal deformations of $(Y,F)$ as a brane are given by pairs $(r, B)\in \Gamma(E^*)\oplus \Omega^2(Y)$ of the form
\begin{align*}
 \label{eq:secondsystemLemma} r&= \rho dq,    \\
 B&=\left(B_{N,0}+d_N\int_0^q\gamma_{N,q'}\right)+dq\wedge \gamma_{N,q}
\end{align*}
where, using the notation $\gamma_{N,q}:=I^*(d_N \rho|_{N_q})$, 
\begin{align}
\rho&\in C^{\infty}(Y)  \text{ satisfies } d_N\int_{0}^1\gamma_{N,q}=0 \text{ for all $q\in S^1$},\\
B_{N,0}&\in \Omega^2(N_0) \text{ is of type $(1,1)$ and closed}. \nonumber
\end{align}
Here $\int_{0}^1\gamma_{N,q}$ denotes the averaging of a family of 1-forms on $N$.
\end{lem}

\begin{proof}
The infinitesimal deformations of $(Y,F)$ as a brane are given by pairs $(r,B)$ as in Corollary \ref{rem:infdefcodim1}.

Write $r=\rho dq$ for $\rho\in C^{\infty}(Y)$.
Write $$B=B_{N,q}+dq\wedge \gamma_{N,q}\in \Omega^2(Y)$$ with $B_{N,q}\in \Omega^2(N_q)$ and $\gamma_{N,q}\in \Omega^1(N_q)$, for all $q\in S^1$.  Condition iii) in Corollary \ref{rem:infdefcodim1} is  equivalent to 
\begin{equation}\label{dq:gammaNq}
    \gamma_{N,q}=I^*(d_N \rho|_{N_q}).
\end{equation}
Hence, in view of Remark \ref{rem:codim1} ii),  $(r,B)$ is an infinitesimal deformation if{f} conditions ii) and iv) in Corollary \ref{rem:infdefcodim1} are satisfied, i.e.\ if{f} 
\begin{itemize}
    \item 
$dB=0$, 
\item $B_{N,q}$ is of type $(1,1)$ for all $q$.
\end{itemize}
Writing the de Rham differential of $Y$  in terms of the factors $N$ and $S^1$, we see that $dB=0$ if{f}  
$$
d_NB_{N,q}=0\quad\quad\text{and}\quad\quad
d_N\gamma_{N,q}=\pd{q}B_{N,q}. 
$$
These equations are satisfied  for all $q\in S^1$ if{f} 
\begin{equation}\label{eq:secondsystem}
d_NB_{N,0}=0\quad\quad\text{and}\quad\quad
d_N\int_{0}^1\gamma_{N,q}=0.
\end{equation}
Indeed, the direct implication holds by Stokes' theorem;  to show the reverse implication, define 
\begin{equation}\label{eq:Bnq}
B_{N,q}=B_{N,0}+d_N\int_0^q\gamma_{N,q'}.
\end{equation}
 Notice that if
$B_{N,0}$ is of type $(1,1)$,
then $B_{N,q}$ is of type $(1,1)$ for all $q\in S^1$, since applying $d_NI^*d_N$ to any function we obtain a form  of  type $(1,1)$, as we saw in Remark \ref{rem:spacefillfunctions}. 
\end{proof}

\begin{proof}[Proof of Proposition \ref{prop:notsurj}]
 In Lemma  \ref{lem:infdefmainex}, the only condition on $r$ is \eqref{eq:secondsystemLemma}, which is just eq. \eqref{eq:imager}.
This shows that if $r$ lies in the image of $\Upsilon$, then it satisfies eq. \eqref{eq:imager}.

Conversely, let $r\in \Gamma(E^*)$ satisfy eq. \eqref{eq:imager}. 
Choose any closed $(1,1)$-form $B_{N,0}$ on $N$ (for instance $B_{N,0}=0$).  Then defining $B$ as in Lemma \ref{lem:infdefmainex} we obtain an infinitesimal brane deformation $(r,B)$, so $r$ lies in the image of $\Upsilon$.
\end{proof}

 A simple instance in which the map $\Upsilon$ is not surjective is the following. Another instance in which $\Upsilon$ fails to be  surjective, in a more extreme way, is given in Example \ref{rmk:K3const} below.
\begin{ex}[$\Upsilon$ is not surjective: $\RR^4$ and $\TT^4$] \label{rem:notsurj}
   We use the fact that any smooth function  $f$ on $N$ can be realized as $\int_{S^1}r$ for some $r\in \Gamma(E^*)$.
As in Example \ref{rem:Idfcoords}, let $N$ be $\RR^4$ or the 4-torus with the 2-forms
given there. Then eq. \eqref{eq:imager} amounts to the system of second order PDE's for $f$ given in that example, and clearly not every smooth function on $N$ satisfies it.  For instance, $f=(x_1)^2$ or $f=(\cos(x_1))^2$ do not satisfy that PDE. 
\end{ex}

\begin{rmk}\label{rem:Lieder}
Eq. \eqref{eq:imager}  is equivalent to $\mathcal{L}_{X_{(  \int_{S^1}r)}}F_N=0$, where $X$ denotes the  hamiltonian vector field w.r.t. $\omega_N$.  This holds by \eqref{eq:LiederXf}.
\end{rmk}

\begin{rmk}
The inclusion ``$\subset$'' in Proposition \ref{prop:notsurj} can be obtained directly differentiating the brane condition given in Proposition \ref{prop:timeone} for a smooth 1-parameter family of branes starting at $(Y,F)$, as we now outline. Let $\{f^t\}\subset  C^{\infty}(Y)$ and
$\{F^t\}\subset  \Omega^2(N)$ be
  smooth 1-parameter families of functions on $Y$ and symplectic forms on $N$, with $f^0=0$, $F^0=F$, so that for every $t$  the  time $1$ flow $\Phi^1_{t}$ of the time-dependent\footnote{Here the role of time is played by the variable $q$.} vector field $\{-X_{{f^t_q}}\}_{q\in [0,1]}$ preserves $F_N^t$.
  Here we use the notation $f^t_q:=f^t|_{N\times\{q\}}$.
Then the foliated 1-form $r:=(\frac{d}{dt}|_0f^t)dq$ satisfies
 equation \eqref{eq:imager}.
 
 To see this,
differentiate at time $t=0$ the equation 
$(\Phi^1_{t})^*F_N^t=F_N^t$ to obtain
$$0=\frac{d}{dt}|_{t=0}(\Phi^1_{t})^*F_N=\cL_{Z_0}F_N,$$ where $\{Z_t\}$ denotes the time-dependent vector field on $N$ corresponding to the 1-parameter family of diffeomorphisms $\{\Phi^1_{t}\}$. One can show  that 
$Z_0$ is the the average over $q\in S^1$ of the vector fields $-\frac{d}{dt}|_{t=0}X_{f^t_q}=-X_{( \frac{d}{dt}|_0f^t_q)}$ on $N$. (To do so, apply \cite[Lemma 3.4]{DoCarmoRG}
to the map $[0,1]^2\to N, (q,t)\mapsto \Phi^q_{t}(x)$ for any point $x\in N$.)
Hence   $Z_0=-X_{(  \int_{S^1}r)}$, and we can conclude the argument using  Remark \ref{rem:Lieder}.
 \end{rmk}

\subsubsection*{Remarks on Proposition 
\ref{prop:notsurj}}

We present several related remarks, and 
Example \ref{rmk:K3const} mentioned above.
The map $\Upsilon$ induces a map in cohomology, whose image forms a subspace of $C^{\infty}(N)$.

\begin{rmk}[Map induced by $\Upsilon$ in cohomology] 
  Since $Y\to N$ is a trivial $S^1$-bundle, the first cohomology $H^1_E$ of the foliated differential forms along $E$ is given by 
the smooth sections of the trivial vector bundle $H^1(S^1)\times N\to N$, i.e.\ by $C^{\infty}(N)$. Under this correspondence, the   class $[r]\in H^1_E$  corresponds to the function $ \int_{S^1}r$. 

 Dividing by   exact degree 1  elements,
  the map $\Upsilon$ induces a map 
  $$\text{First cohomology of the complex \eqref{eq:defcplxinv}}\to  C^{\infty}(N)$$
  whose image is $\{f\in C^{\infty}(N): d_N I^* d_N f=0\}$. 
  An element $r\in \Gamma(E^*)$ lies in the image of $\Upsilon$ if{f} $\int_{S^1}r$ lies in the image of the above map. 
\end{rmk}

By the following two remarks, the  infinitesimal hamiltonian  (respectively    symplectic) deformations  of $Y$ as a coisotropic submanifold induce in cohomology the zero function (resp. constant functions) on $N$.

\begin{rmk}[Hamiltonian deformations]
\label{rem:trivialinfdefham}
Infinitesimal \emph{hamiltonian} deformations
of $Y$ as a coisotropic submanifold
are given by exact elements of
$\Gamma_{\text{cl}}(E^*)$, i.e.\ those of the form $r=d_Eg$ for some $g\in C^{\infty}(Y)$ (see \cite[Prop. 3.5 and Rem. 3.6]{Eqcoiso}). An element $r\in \Gamma_{\text{cl}}(E^*)$ is exact if{f} $\int_{S^1}r=0$. In particular, such an element 
$r$ must lie in the image of $\Upsilon$, by Proposition \ref{prop:notsurj}. This can be seen in a geometric way too: such an element gives rise to a hamiltonian flow on $M$, which one can use to obtain a path of branes in $M$ as in Remark \ref{rmk:branefct}, and thus an  infinitesimal brane deformation that maps to $r$ under $\Upsilon$.
\end{rmk}

\begin{rmk}[Deformations by symplectomorphisms]
\label{rem:trivialinfdef}
Infinitesimal deformations  of $Y$ (as a coisotropic submanifold)  by \emph{symplectomorphisms} of the ambient symplectic manifold 
are given by  elements  
$r\in \Gamma_{\text{cl}}(E^*)$ which are the pullback to the leaves of a closed 1-form $R\in \Omega^1_{\text{cl}}(Y)$ (see \cite[Prop. 4.5 and Rem. 4.6]{Eqcoiso}). 
We claim that such elements $r$ are precisely the elements of $\Gamma_{\text{cl}}(E^*)$ such that $\int_{S^1}r $ is a \emph{constant} function on $N$.  In particular, such an element 
$r$ must lie in the image of $\Upsilon$, by Proposition \ref{prop:notsurj}. This can be seen also using the same geometric argument as in Remark \ref{rem:trivialinfdefham}.   

We now prove the above claim. If $r$ is the restriction of $R\in \Omega^1_{\text{cl}}(Y)$ to leaves, by Stokes' theorem we  have  $$\int_{S^1\times\{x\}}r- \int_{S^1\times\{y\}}r=\int_CdR=0$$
 for any $x,y\in N$, where $C$ is the product of $S^1$ and an arc in $N$ joining $x$ to $y$; therefore, $\int_{S^1}r $ is a  constant  function. Conversely, let be given $r=\rho dq$ for $\rho\in C^{\infty}(Y)$. Any extension can be written as $$R=\rho dq+\alpha_{N,q}\in \Omega^1(Y),$$ where $\alpha_{N,q}\in \Omega^1(N)$ depends smoothly on $q\in S^1$;
 the extension is closed if{f}
$d_N\alpha_{N,q}=0$ and $d_N \rho=\pd{q}\alpha_{N,q}$. 
Now assume that $\int_{S^1}r $ is constant, and set   $\alpha_{N,q}:=d_N (\int_0^q \rho|_{N\times \{q'\}})$. This gives a well-defined\footnote{Indeed $\alpha_{N,0}=\alpha_{N,1}$, as they both vanish.} 1-form $R$ on $Y$,  which furthermore is closed.
\end{rmk}

When $N$ is a $K3$ manifold, all infinitesimal brane deformations are infinitesimal deformations by symplectomorphisms, as we now show.

\begin{ex}[$\Upsilon$ is not surjective: the K3 manifold]
    \label{rmk:K3const}
Assume that $N$ is the K3 manifold.   Then any infinitesimal deformation $(r,B)$
of $Y$ as a brane
is such that the function 
$\int_{S^1}r$ on $N$ is constant, as we explain below.  By Remark \ref{rem:trivialinfdef}, it is   therefore an infinitesimal deformation of $Y$ by symplectomorphisms.
We do not know whether
all branes nearby $Y$ are obtained from $(Y,F)$ applying a symplectomorphism of $M$, see Remark \ref{rmk:K3sympl}. 

We now prove the above statement. By Proposition \ref{prop:notsurj} and Remark \ref{rem:Lieder} we know that the vector field $X_{(  \int_{S^1}r)}$ on $N$ preserves $F_N$. Being a hamiltonian vector field, it preserves $\omega_N$ too, hence it preserves $I$. By Proposition \ref{prop:K3id} in \S \ref{subsec:notbrane} we know that it must be the zero vector field, i.e.\ $\int_{S^1}r$ is a constant function.
\end{ex}

 
\appendix

\section{Two proofs}
\label{app:proof}

This appendix contains the proofs of {three} lemmas in the body of the paper. 

\subsection{Proof of Lemma \ref{rem:Melanie}}

Lemma \ref{rem:Melanie} follows immediately from the following more general lemma, in which the distribution $G$ is not assumed to be involutive.

\begin{lem}\label{rem:Melanie0}
 Let $Y$ be a manifold and $\omega$ a 2-form of constant rank, denote $E:=\ker(\omega)$. Let $G$ be a distribution such that  $E\oplus G=TY$.
Then $\omega$ is closed if{f}:
\begin{itemize}
\item[i)] the distribution $E$ is involutive
\item[ii)]  $(\cL_X\omega)|_{\wedge^2 G}=0$ for all $X\in \Gamma(E)$
\item[iii)] $(d\omega)|_{\wedge^3 G}=0$.
\end{itemize}
\end{lem}
\begin{proof}  
We split $\wedge^3 TY$ as a direct sum  $\wedge^3 E\oplus (\wedge^2 E\otimes G)\oplus  (E\otimes \wedge^2G)\oplus\wedge^3 G$,
and use the well-known formula in terms of Lie brackets to compute $d\omega$. We see that $ (d\omega)|_{\wedge^3 E}$ automatically vanishes. Assume that $d\omega=0$. Then:

i) follows because for all $X,Y\in \Gamma(E)$ and $Z^{\perp}\in \Gamma(G)$ we have $d\omega(X,Y, Z^{\perp})=
-\omega([X,Y], Z^{\perp})$.

 ii) follows because for all $X\in \Gamma(E)$ and $Y^{\perp}, Z^{\perp}\in \Gamma(G)$ we have $$d\omega(X,Y^{\perp}, Z^{\perp})=\cL_X\left(\omega(Y^{\perp}, Z^{\perp})\right)-\omega([X,Y^{\perp}], Z^{\perp})
+\omega([X,Z^{\perp}], Y^{\perp})=(\cL_X\omega)(Y^{\perp}, Z^{\perp}).$$
Here we use the Leibniz rule for the Lie derivative in the second equality.  

iii) follows trivially.
\end{proof}

\subsection{Proof of Lemma \ref{lem:infdefdirect}}

The proof of Lemma \ref{lem:infdefdirect} is based on ideas by Collier \cite[Prop. 19, page 46]{Collier}, which rely on  generalized complex geometry.
We therefore first recall the notion of brane in that setting.
 
Let $\mathcal{J}$ be a  generalized complex structure  on the standard   Courant algebroid $TM\oplus T^*M$.
A \emph{brane} $(Y,F)$ consists of a submanifold $Y$ endowed with a closed 2-form $F$
such that the generalized tangent bundle
\begin{equation}\label{eq:tauF}
    \tau_F^Y := \{ (X,\xi) \vert X \in TY, \xi \in T^*M|_Y \text{ s.t. } \xi|_{TY}=\iota_X F\}  
\end{equation}  
satisfies $\mathcal{J}(\tau_F^Y)=\tau_F^Y$, see \cite{GualtieriThesis}. 
Notice that $\tau_F^Y$ is lagrangian w.r.t. the symmetric pairing $\langle \cdot,\cdot \rangle$ on $TM\oplus T^*M$.
For generalized complex structures that arises from a symplectic form $\widetilde{\omega}$ on $M$, 
i.e.\ $ \mathcal{J}=\left( \begin{matrix} 
0 & -\widetilde{\omega}^{-1} \\ 
\widetilde{\omega} & 0 
\end{matrix}\right)$, branes can be characterized as in Definition \ref{def:brane}. 
We will use the  notation $$e^B(X,\xi):=(X,\xi+\iota_X B)$$ for $B\in \Omega^2(M)$ and $(X,\xi)\in TM\oplus T^*M$.

\begin{proof}[Proof of Lemma \ref{lem:infdefdirect}]
We use freely the notation introduced just before Lemma \ref{lem:infdefdirect}. In particular we assume that $M$ is a neighborhood of the zero section in the  total space of the  vector bundle $q\colon E^*\to Y$. Further we denote by $\tilde{\omega}$ the ambient symplectic form and by $\omega$ its pullback to $Y$.

We know that $(Y_t,F_{Y_t})$ is a curve of brane deformations if{f}   
 $\mathcal{J}$ preserves $\tau^{Y_t}_{F_{Y_t}}$ for all $t$. Let $\phi_t\colon E^*\to E^*$ be the translation by $-s_t$; it maps $Y_t$ to the zero section $Y$ of $E^*$ and $F_{Y_t}$ to $F_t$, so it maps $\tau^{Y_t}_{F_{Y_t}}$ to $\tau^{Y}_{F_t}$.
 Hence $(\phi_t)_*\mathcal{J}$ preserves $\tau^{Y}_{F_t}:=\tau_{F_t}$. Since
the latter is a lagrangian subbundle (over $Y$), this can be expressed as follows: for all $\mathbb{X}=(X,\xi), \mathbb{Z}=(Z,\eta)\in \tau_F$, and for all $t$, 
\begin{equation}\label{eq:branecondDirac}
0=    \left\langle ((\phi_t)_*\mathcal{J})(e^{\widetilde{F}_t-\widetilde{F}}\mathbb{X}), e^{\widetilde{F}_t-\widetilde{F}} \mathbb{Z} \right\rangle.
\end{equation}
Here we fix extensions $\widetilde{F_t}\in \Omega^2(E^*)$ of $F_t$  to the ambient space, allowing to write $\tau_{F_t}=e^{\widetilde{F}_t-\widetilde{F}}\tau_{F}$.

Taking the time derivative at zero  of eq. \eqref{eq:branecondDirac} we get 
\begin{equation}\label{eq:ddtCollier}
0=  \left\langle  \left(\frac{d}{dt}|_0   (\phi_t)_*\mathcal{J} \right)\mathbb{X},\mathbb{Z}\right\rangle+
 \left\langle \mathcal{J}\left(\frac{d}{dt}|_0  e^{\widetilde{F}_t-\widetilde{F}}\mathbb{X}\right),   \mathbb{Z} \right\rangle+
 \left\langle  \mathcal{J}( \mathbb{X}), \frac{d}{dt}|_0  e^{\widetilde{F}_t-\widetilde{F}} \mathbb{Z} \right\rangle. 
\end{equation}
Now  $ (\phi_t)_*\mathcal{J}=\left( \begin{matrix} 
0 & -\widetilde{\omega_t}^{-1} \\ 
\widetilde{\omega_t} & 0 
\end{matrix}\right)$
where $\widetilde{\omega}_t$ is the pullback  by $(\phi_t)^{-1}$ of  the ambient symplectic form $\widetilde{\omega}$ on $E^*$; therefore 
$\frac{d}{dt}|_0   (\phi_t)_*\mathcal{J}=\left( \begin{matrix} 
0 &  \widetilde{\omega}^{-1}\dot{\widetilde{\omega}}{\widetilde{\omega}}^{-1} \\ 
\dot{\widetilde{\omega} }& 0 
\end{matrix}\right)$. 
Further $\frac{d}{dt}|_0  e^{\widetilde{F}_t-\widetilde{F}}\mathbb{X}=(0,\dot{\widetilde{F}}X)$.
Hence equation \eqref{eq:ddtCollier} reads
 \begin{align*}
0&=   \left\langle \left( \begin{matrix} 
 \widetilde{\omega}^{-1}\dot{\widetilde{\omega}}{\widetilde{\omega}}^{-1} \xi \\ 
  \dot{\widetilde{\omega}} X
\end{matrix}\right),
\mathbb{Z}
 \right\rangle
 +
    \left\langle \left( \begin{matrix} 
 -\widetilde{\omega}^{-1} \dot{\widetilde{F}}X \\ 
 0
\end{matrix}\right),
\mathbb{Z}
 \right\rangle
 +
     \left\langle \left( \begin{matrix} 
  -\widetilde{\omega}^{-1} \xi\\ 
 \widetilde{\omega}X
\end{matrix}\right)
,
  \left( \begin{matrix} 
  0\\ 
\dot{\widetilde{F}}Z
\end{matrix}\right) \right\rangle\\
     &= \langle   \widetilde{\omega}^{-1}\dot{\widetilde{\omega}}{\widetilde{\omega}}^{-1} \xi\;,\;  \eta\rangle+
 \langle \dot{\widetilde{\omega}} X \;,\;Z \rangle
-   \langle   \widetilde{\omega}^{-1} \dot{\widetilde{F}}X \;,\;    \eta \rangle
+     \langle   \dot{\widetilde{F}}\widetilde{\omega}^{-1} \xi \;,\;  Z \rangle.
   \end{align*}
Since $\tau_F$ is lagrangian, this is equivalent to the following:
\begin{equation}\label{eq:defbranecond}
  \left( \begin{matrix} 
  \widetilde{\omega}^{-1}(\dot{\widetilde{\omega}}{\widetilde{\omega}}^{-1} \xi-
  \dot{\widetilde{F}}X)\\ 
 \dot{\widetilde{\omega}} X+  \dot{\widetilde{F}}\widetilde{\omega}^{-1} \xi
\end{matrix}\right) 
\in \tau_F \quad\quad\quad\text{for all $(X,\xi)\in \tau_F$}.
\end{equation}

We now show how eq. \eqref{eq:defbranecond} corresponds to the items in the statement. We will use that,  by the coisotropicity of $Y$, the isomorphism $\widetilde{\omega}$ maps $TY$ to $E_{amb}^{\circ}$
(the annihilator of the distribution $E$ in the ambient manifold $E^*$)
and $E$ to $TY^{\circ}$. 

\emph{Item i)} and \emph{item ii)}: 
 by taking $X=0$, in view of \eqref{eq:tauF} the covector $\xi$ ranges over all elements of  $TY^{\circ}$, so that $\widetilde{\omega}^{-1} \xi$ ranges over all elements of $E$. In this case,
the tangent component of eq. \eqref{eq:defbranecond} lying in $TY$ implies that $\dot{\omega}\in \Omega^2(Y)$ is a horizontal form, i.e.\ vanishes on $\wedge^2 E$. Further, in this case, the condition on the cotangent component in eq. \eqref{eq:defbranecond}  implies that $\dot{F}\in \Omega^2(Y)$ is also a  horizontal form, using the fact that $E=\ker(F)$. 
(The fact that both forms are horizontal also follows from  \cite{DefPres}.)

 We saw in \S \ref{subsec:Diracapproach} that the     $Y_t$ being coisotropic submanifolds implies \eqref{eq:OPomegat}, from which we obtain   $\dot{\omega}=-d\bar{\dot{s}}$. 
Beside the point-wise conditions, we also have to take into account the fact
the $F_t$  are  closed forms,
implying that $\dot{F}$ is closed. Thus we obtain items \emph{i)} and \emph{ii)}.

In the following we will use repeatedly the following observations:
\begin{itemize}
    \item[$(\star)$]
For any $(X,\xi)\in \tau_F$, 
the projection  $TY\to TY/E$ maps  $\widetilde{\omega}^{-1} \xi\in TY$  to $I[X]$. 

  \item[$(\star\star)$] If $\xi \in E_{amb}^{\circ}$ then $F  (\widetilde{\omega}^{-1}\xi)=I^*(\xi|_{TY})$.
\end{itemize}

 (To show the first one, since $\omega\circ I= F$ on $TY/E$, we have to check that $\omega(\widetilde{\omega}^{-1} \xi)=FX$, which holds since $\omega(\widetilde{\omega}^{-1}\xi)=
( \widetilde{\omega}(\widetilde{\omega}^{-1}\xi))|_{TY}=\xi|_{TY}=\iota_XF$. For the second one, we have $(\xi|_{TY}) \in E^{\circ}$ and $I^*(\xi|_{TY})=(F\circ \omega^{-1})(\xi|_{TY})$, and use that a lift of $\omega^{-1}(\xi|_{TY})$ to $TY$ is  $\widetilde{\omega}^{-1}\xi$.)

\emph{Item iii)}: 
The condition on the  tangent component of eq. \eqref{eq:defbranecond} is that it lies in $TY$. Upon applying the isomorphism $\widetilde{\omega}$,  this condition can be equivalently written as: 
\begin{equation}\label{eq:tgEcirc}
\dot{\widetilde{\omega}}{\widetilde{\omega}}^{-1} \xi-  \dot{\widetilde{F}}X \in E_{amb}^{\circ}.
\end{equation} 
In turn this means that $
 \dot{\omega}(\widetilde{\omega}^{-1}\xi,e)=\dot{F}(X,e)$ for all $e\in E$, which by   observation $(\star)$ is exactly eq. \eqref{eq:samekernelDirac}.

\emph{Item iv)}: The condition on the cotangent component of eq. \eqref{eq:defbranecond} reads
\begin{equation}\label{eq:ctgtcomp}
(\dot{\widetilde{\omega}} X+  \dot{\widetilde{F}}\widetilde{\omega}^{-1} \xi)|_{TY}=
F( \widetilde{\omega}^{-1}(\dot{\widetilde{\omega}}{\widetilde{\omega}}^{-1} \xi- \dot{\widetilde{F}}X)).
\end{equation}
The l.h.s. lies in $E^{\circ}$, by  \eqref{eq:defbranecond} and the fact that whenever $(X',\xi')\in \tau_F$, then $\xi' \in E_{amb}^{\circ}$, since $E=\ker(F)$. 
The r.h.s. of the above equation 
 can be written as 
 $I^* ((\dot{\widetilde{\omega}}{\widetilde{\omega}}^{-1} \xi-  \dot{\widetilde{F}}X)|_{TY})$,
 thanks to  eq. \eqref{eq:tgEcirc} and observation $(\star\star)$. 

Hence  
\eqref{eq:ctgtcomp}
can be written as the following equation in $E^{\circ}\subset T^*Y$:
\begin{equation}
  \dot{\omega} X+  \dot{F}\widetilde{\omega}^{-1} \xi =I^*(\dot{\omega}(\widetilde{\omega}^{-1} \xi)-\dot{F}X).
\end{equation}
This is exactly eq. \eqref{eq:I^2Dirac}, as can be seen using again observation $(\star)$.
\end{proof}

\subsection{Proof of Lemma \ref{lem:T4unobstr}}

Lemma \ref{lem:T4unobstr} is expected since a certain moduli space of space-filling brane structures on the $4$-torus is smooth  \cite{KLZTorelli}. Here we present a proof based on more elementary arguments from \cite{KLZTorelli}. 

\begin{proof}[Proof of Lemma \ref{lem:T4unobstr}]
We may assume that  $\omega$, $F$ and $I$ are all constant on the torus, by the first claim in the proof of \cite[Prop. 4.17]{KLZTorelli}.

Consider the codimension $2$ submanifold   $\mathcal{Q}_{[\omega]}$  of $H^2(M,\RR)$ consisting of elements $\phi$ which satisfy
    \begin{equation}\label{eq:QR}
            \phi\wedge [\omega]= 0 ,\quad \phi \wedge \phi=[\omega] \wedge [\omega].  
    \end{equation}  
The tangent space to $\mathcal{Q}_{[\omega]}$ at $[F]$ is $H_{\RR}^{1,1}(M)$, by \cite[Rem. 4.16]{KLZTorelli}.

Now let $B$ be an infinitesimal deformation of $F$, i.e.\ $B\in \Omega^{(1,1)}_{\RR,\text{cl}}$  (see Lemma \ref{lem:infdefspace-filling}).
Take a smooth curve $\phi(t)$ in $\mathcal{Q}_{[\omega]}$ through $[F]$ such that $\frac{d}{dt}|_0 \phi(t)=[B]$. Let $F(t)\in \Omega^2(M)$ be the unique constant representative of the class $\phi(t)$. Then each $F(t)$ is a space-filling brane structure by \cite[Prop. 4.1]{KLZTorelli}, since it satisfies  $F(t)\wedge \omega= 0$ and  $F(t) \wedge F(t)=\omega \wedge \omega$. Further, $\frac{d}{dt}|_0F(t)=\widehat{B}$, where
$\widehat{B}  \in \Omega^{1,1}_{\RR}(M)$ is the unique constant representative of the class $[B]\in H^{1,1}(M)_{\RR}\subset H^2(M)$. Hence the  infinitesimal deformation 
$\widehat{B}$ can be prolonged. By Remark \ref{rem:additionalnote} and Lemma \ref{lem:equivprolong}, the same holds for $B$.
\end{proof}


 \bibliographystyle{halpha} 

\begin{thebibliography}{BPVdV84}

\bibitem[BPVdV84]{BarthPetersVandeVen1984}
W.~Barth, C.~Peters, and A.~Van~de Ven.
\newblock {\em Compact complex surfaces}, volume~4 of {\em Ergebnisse der
  Mathematik und ihrer Grenzgebiete (3) [Results in Mathematics and Related
  Areas (3)]}.
\newblock Springer-Verlag, Berlin, 1984.

\bibitem[CFM21]{PoisGeoBookAMS}
Marius Crainic, Rui~Loja Fernandes, and Ioan M\u{a}rcu\c{t}.
\newblock {\em Lectures on {P}oisson geometry}, volume 217 of {\em Graduate
  Studies in Mathematics}.
\newblock American Mathematical Society, Providence, RI, [2021] \copyright
  2021.

\bibitem[Col14]{Collier}
Braxton~L. Collier.
\newblock Deformations of generalized complex branes.
\newblock 03 2014, Preprint ArXiv:1403.2970.

\bibitem[Cou90]{Cou}
Ted Courant.
\newblock Dirac manifolds.
\newblock {\em Trans. Amer. Math. Soc.}, 319(2):631--661, 1990.

\bibitem[dC92]{DoCarmoRG}
Manfredo Perdig\~{a}o do~Carmo.
\newblock {\em Riemannian geometry}.
\newblock Mathematics: Theory \& Applications. Birkh\"{a}user Boston, Inc.,
  Boston, MA, portuguese edition, 1992.

\bibitem[Deb00]{DebordThesis}
Claire Debord.
\newblock Feuilletages singuliers et groupo\"ides d'holonomie.
\newblock {Phd Thesis at Institut de Math\'ematiques de Jussieu}, 2000.

\bibitem[Got82]{Gotay}
Mark~J. Gotay.
\newblock On coisotropic imbeddings of presymplectic manifolds.
\newblock {\em Proc. Amer. Math. Soc.}, 84(1):111--114, 1982.

\bibitem[Gua03]{GualtieriThesis}
Marco Gualtieri.
\newblock {\em Generalized complex geometry}.
\newblock PhD thesis, University of Oxford, 2003, Preprint ArXiv:math/0401221.

\bibitem[GW09]{GukovWitten}
Sergei Gukov and Edward Witten.
\newblock Branes and quantization.
\newblock {\em Adv. Theor. Math. Phys.}, 13(5):1445--1518, 2009.

\bibitem[Huy05]{Huy}
Daniel Huybrechts.
\newblock {\em Complex geometry}.
\newblock Universitext. Springer-Verlag, Berlin, 2005.

\bibitem[{Huy}16]{HuyLecturesK3}
Daniel {Huybrechts}.
\newblock {\em {Lectures on \(K\)3 surfaces}}, volume 158.
\newblock Cambridge University Press, 2016.

\bibitem[KLZ25]{KLZTorelli}
Charlotte Kirchhoff-Lukat and Marco Zambon.
\newblock Moduli spaces of spacefilling branes in symplectic 4-manifolds.
\newblock {\em Mathematische Zeitschrift}, 311(1):15, 2025.

\bibitem[KM06]{KoerMar}
Paul Koerber and Luca Martucci.
\newblock {Deformations of calibrated D-branes in flux generalized complex
  manifolds}.
\newblock {\em Journal of High Energy Physics}, 2006.

\bibitem[KO03]{KO}
Anton Kapustin and Dmitri Orlov.
\newblock Remarks on {A}-branes, mirror symmetry, and the {F}ukaya category.
\newblock {\em J. Geom. Phys.}, 48(1):84--99, 2003.

\bibitem[LWX97]{LWX}
Zhang-Ju Liu, Alan Weinstein, and Ping Xu.
\newblock Manin triples for {L}ie bialgebroids.
\newblock {\em J. Differential Geom.}, 45(3):547--574, 1997.

\bibitem[OP05]{OP}
Yong-Geun Oh and Jae-Suk Park.
\newblock Deformations of coisotropic submanifolds and strong homotopy {L}ie
  algebroids.
\newblock {\em Invent. Math.}, 161(2):287--360, 2005.

\bibitem[SZ17]{Eqcoiso}
Florian Sch{\"a}tz and Marco Zambon.
\newblock Equivalences of coisotropic submanifolds.
\newblock {\em J. Symplectic Geom.}, 15(1):107--149, 2017.

\bibitem[SZ20]{DefPres}
Florian Sch\"{a}tz and Marco Zambon.
\newblock Deformations of pre-symplectic structures and the {K}oszul
  {$L_\infty$}-algebra.
\newblock {\em Int. Math. Res. Not. IMRN}, (14):4191--4237, 2020.

\end{thebibliography}

\end{document}